\documentclass[12pt]{article} 

\usepackage{amsfonts}
\usepackage{amsmath}
\usepackage{amssymb}
\usepackage{graphicx}

 \input epsf

\newtheorem{theorem}{Theorem}

\newtheorem{proposition}[theorem]{Proposition}

\newtheorem{corollary}[theorem]{Corollary}

\newtheorem{lemma}[theorem]{Lemma}

\newtheorem{remark}[theorem]{Remark}

\newtheorem{example}[theorem]{Example}

\newtheorem{definition}[theorem]{Definition}

\def\CC{{\mathbb C}}
\def\NN{{\mathbb N}}
\def\RR{{\mathbb R}}

\def\qed{\hfill \vrule height7pt width3pt depth0pt}

\def\ux{\underline x}

\def\bfe{{\bf e}}
\def\bfL{{\bf L}}
\def\bfH{{\bf H}}

\def\dive{{\rm div}\,}
\def\curl{{\rm curl}\,}

\numberwithin{equation}{section}
\numberwithin{theorem}{section}

\begin{document}

\title{Optimal Focusing for Monochromatic Scalar and
Electromagnetic Waves}

\author{\sc \small  Jeffrey Rauch\thanks{University of Michigan, Ann Arbor 48109 MI, USA. email: rauch@umich.edu.  Research partially supported by NSF under grant NSF DMS-0807600} 
}

\date{}

\maketitle

\begin{abstract}For monochromatic solutions of D'Alembert's
wave equation and Maxwell's equations, we obtain
sharp bounds on the sup norm as a function of the 
far field energy.  The extremizer in the scalar case
is radial.  In the case of Maxwell's equation, the electric
field maximizing the value at the origin
follows longitude
lines on the sphere at infinity.  
In dimension $d=3$ the highest electric field 
for Maxwell's equation is smaller by a factor 2/3
than the highest corresponding scalar waves.

The highest electric field densities on the balls  $B_R(0)$ occur
as $R\to 0$.  The density dips to half max at $R$
approximately equal to one third the wavelength.
The extremizing fields are 
identical to those
that attain the maximum field intensity at the 
origin.
\end{abstract}



\section{Introduction.}

The problem we address is to find for fixed frequency $\omega/2\pi$,
the 
monochromatic
solutions of the wave equation and of Maxwell's
equation which achieve the highest field values 
at a point or more generally the greatest electrical energy in
 ball
of fixed small radius.   They are constrained by the 
energy at $|x|\sim\infty$.
This leads to several
variational problems;

\vskip.1cm

$\bullet$ Maximize
the field strength at a  point.

$\bullet$ For fixed $R$, maximize the energy 
in a ball of radius $R$.

\vskip.1cm

A third problem is,

\vskip.1cm

$\bullet$  Find $R_{optimal}$ so that the
energy {\sl density} is largest.

\vskip.1cm
We show that the last is degenerate, the maximum
occuring at $R=0$.

There are experiments in course whose strategy
is to focus a number of  coherent
high power laser beams on a small volume
to achieve very high energy densities.  The problem was proposed to me by
G. Mourou because of his leadership role in  the 
European Extreme
Light Initiative.  If by better focusing one
can reduce the size of the incoming lasers there would
be significant
benefits.  Identification of the extrema 
guides the deployment of the lasers.
With the experimental design as motivation the 
maximization for the Maxwell equations is
interpreted as maximization of the energy in the 
electric field $E$, ignoring the magnetic
contribution.  Including the magnetic contribution
creates an analogous problem amenable to
the techniques introduced here.

As natural as these questions appear, I have been
unable to find  previous work on them.

Study solutions of the
scalar wave equation and  of Maxwell equations,
$$
v_{tt}\ -\ \Delta\,v \ =\ 0\,,
\qquad
E_t = \curl B,
\quad
B_t=-\curl E\,,
\quad
\dive E =\dive B = 0\,,
$$
for spatial dimensions $d\ge 2$.  Units are chosen
so that the propagation speed is 1.

\begin{definition}
A solution $v$ of the wave equation is monochromatic
if it is of the form
\begin{equation}
\label{eq:mono}
v\ =\ 
\psi(t) \, 
u(x),
\qquad 
{\rm with}
\qquad
\psi^{\prime\prime}  +\ 
\omega^2\,\psi \ =\ 0\,,
\quad
\omega\
>\
0\,.
\end{equation}
Monochromatic  solutions 
of Maxwell's equation are those of the form
\begin{equation}
\label{eq:mono2}
\psi(t) \, 
\big(
E(x)\,,\,B(x)
\big)
\,,
\qquad
{\rm with }
\qquad
\psi^{\prime\prime}  +\ 
\omega^2\,\psi \ =\ 0\,,
\quad
\omega\
> \
0\,.
\end{equation}
\end{definition}

They are generated by
\begin{equation}
\label{eq:monoplusminus}
e^{\pm i\omega t} \, u(x),
\qquad 
{\rm and}
\qquad
e^{\pm i\omega t}\,\big(E(x)\,,\,B(x)\big)\,.
\end{equation}
Scaling $t,x\to \omega t,\omega x$ reduces the study to the case
$\omega=1$.
In that case, the {\bf reduced wave equations} are satisfied,
\begin{equation}\label{eq:reduced}
(\Delta +1)u(x)\ =\ 0\,,
\qquad
(\Delta +1)E(x)
\ =\ 
(\Delta +1)B(x)\ =\ 0\,.
\end{equation}

\vskip.2cm

{\bf Notation.}  {\sl  The absolute value sign $|\ |$ is used to denote
the modulus of complex numbers, the length of vectors in $\CC^d$,
surface area, and, volume. 
}
Examples:  $|S^{d-1}|$ and $|B_R(0)|$.

\begin{example}
\label{ex:planewaves}
   The plane waves $e^{i(\pm\omega t  +  \xi x)}$
with $|\xi|=\omega$ is a monochromatic solution of the wave
equation.  Its period and wavelength are equal to
 $2\pi/\omega$.  For Maxwell's equations the analogue
 is $E=e^{i(\pm\omega t  +  \xi x)}\bfe$ with $\bfe\in \CC^d$
 satisfying $\xi\cdot\bfe=0$ to guarantee the divergence free
 condition. 
 \end{example}

The solutions that interest us tend to zero as $|x|\to \infty$.

\begin{example}  When   $d=3$, 
$
u(x):=
\sin  |x|/|x|
$
is a solution of the reduced wave equation
(see also Example \ref{ex:scalmax}).  The corresponding
solutions of the wave equation is
$$
v\ =\ 
e^{\pm it} 
\frac{\sin |x|}
{|x|}
\ =\ 
\frac{1}{2}
\Big(
\frac
{e^{i(\pm t+|x|) } }
{|x|}
\ -\
\frac
{ e^{ i(\pm t-|x|) } }
{|x|}
\Big)\,.
$$
For the plus sign, 
the first term represents an incoming spherical wave
and the second 
outgoing.  To create such a solution it suffices
to generate the incoming wave.   The outgoing wave with the change
of sign is then generated by that wave after it focuses at the
origin.
\end{example}

\begin{example}
Finite energy solutions of Maxwell's equations
are those for which $\int_{\RR^d} |E|^2 + |B|^2\,dx<\infty.$
They 
satisfy
$$
\forall R>0,\qquad
\lim_{t\to \infty}
\ 
  \int_{|x|\le R} |E(t,x)|^2 + |B(t,x)|^2\ dx
\ =\
0\,.
$$ 
Therefore,
 the solution $(E,B)=0$ is the 
only monochromatic solution of finite energy.
\end{example}
 
The solutions $(E(x),B(x))$ that
 tend
to zero as $x\to \infty$ define
tempered distributions on $\RR^d$.
When $(E(x),B(x))$ 
is a tempered solution of 
the reduced wave equation,
the Fourier Transforms
satisfy
$$
(1-|\xi|^2)\widehat E(\xi)
\ =\ 
(1-|\xi|^2) \widehat B(\xi) \ =\ 0.
$$
Therefore the support of $\widehat E$ is contained in the unit
sphere $S^{d-1}:=\big\{|\xi|=1\big\}$.  Since $1-|\xi|^2$
has nonvanishing gradient on this set it follows that
the value of $\widehat E$ on a test function $\psi(\xi)$ is 
determined by the restriction of $\psi$ to $S^{d-1}$.
Therefore there is a distribution $\bfe\in {\cal D^\prime}(\{|\xi|=1\})$
so that
\begin{equation}
\label{eq:E}
E(x) 
\ :=\ 
\int_{|\xi|=1}
e^{ix\xi}
\
\bfe(\xi)\
d\sigma\,,
\end{equation}
where we use the usual abuse of notation indicating as an
integral the pairing of the distribution $\bfe$ with the 
test function $e^{ix\xi}\big|_{|\xi|=1}$.
Conversely, every such expression
is a tempered vector valued solution of the
 reduced wave equation.  

If $\bfe$ is smooth, the principal of stationary phase
(see \S \ref{sec:competing})
shows that 
as
$|x|\to \infty$,
\begin{equation}
\label{eq:stationary}
E(x)  =
\frac{1/\sqrt{2\pi}}
{
|x|^{(d-1)/2}
 }
\
\Big(
e^{-i | x |}\,\bfe((-x/|x|)
+
e^{i\pi(d-1)/4}\, e^{i|x|}\,
\bfe(x/|x|)
+ O(1/|x|)\Big)\,.
\end{equation}
The field in $O(|r|^{-(d-1)/2})$.  In particular,
\begin{equation}
\label{eq:slowgrow}
\sup_{R\ge 1}\ 
R^{-1}\int_{|x|\le  R}|E(x)|^2\ d\sigma
\ <\
\infty\,,\quad
{\rm and.}
\end{equation}
\begin{equation}
\label{eq:L2}
\lim_{R\to \infty}
\int_{R\le |x|\le 2R}
|E(x)|^2\ dx
\
=\
c_d\, \int_{|\xi|=1} |\bfe(\xi)|^2\ d\sigma\,.
\end{equation}
For a field defined by a distribution $\bfe$,
\eqref{eq:slowgrow}  holds
if and only if
$\bfe\in L^2( S^{d-1} )$.
  In 
  that case \eqref{eq:L2} holds and 
 stationary phase approximation holds in 
 an $L^2$ sense.
 This is the class of solutions of the reduced wave
equation that we study.  Equation \eqref{eq:L2}
shows that 
$ \| \bfe \|_{L^2(S^{d-1})} $ 
is a natural measure
of the strength of the field at infinity.

The divergence free condition 
in  Maxwell's
equations is satisfied
if and only if $\xi\cdot \bfe(\xi)=0$ on $S^{d-1}$.
In that case,
the solutions 
$e^{\pm it}E(x)$ of the 
 time dependent equation   are 
 linear combinations of the plane waves in 
 Example \ref{ex:planewaves}.
Denote by $\bfH$ the closed subspace of 
$\bfe\in L^2(S^{d-1};\CC^d)$ with $\xi\cdot\bfe=0$.

For $\bfe\in \bfH$ and $x=r\xi$ with $|\xi|=1$ and $r>>1$,
the solution $e^{it}E(x)$ of Maxwell's 
equations
satisfies
$$
r^{ (d-1)/2 } \, E(x)
\ \approx\
\frac{1}
{\sqrt{2\pi} }
\Big(
e^{ i(t-r) } \,\bfe(-\xi)
\ +\ 
e^{ i\pi(d-1)/4 }\,e^{ i(t+r) } \, \bfe(\xi)
\Big)\,.
$$
In practice, the incoming wave
$$
\frac{1}
{ \sqrt{2\pi} }
\  
\int_{ S^{d-1} }
e^{ i\pi(d-1)/4 }\
\frac{ e^{ i(t+r) } }
{ r^{ (d-1)/2 } } \ \bfe(\xi)\ d\sigma
$$
is generated at large $r$ and the monochromatic solution is
observed for $t>> 1$.  The phase factor $e^{i\pi(d-1)/4 }$
correponds to  the 
phase shift from the focusing at the origin.

The first two variational problems 
for Maxwell's equations
seek to maximize
$$
J_1(\bfe) \ :=\
|E(0)|^2\,,
\qquad
{\rm and}
\qquad
J_2(\bfe) \ :=\ 
\int_{|x|\le R} |E(x)|^2\ dx\,.
$$
among $\bfe\in \bfH$ with $\int_{S^{d-1}}|\bfe(\xi)|^2\,d\sigma=1$.

Theorems \ref{thm:scalarmaxpoint} and \ref{thm:maxpoint} compute
the maxima of $J_1$
in the scalar and electromagnetic cases.
The maximum in the scalar case and also the vector case
wiithout divergence free condition is $|S^{d-1}|$.
It is attained when and only when $\bfe$ is constant. 
 For the electromagnetic
case, $\xi\cdot \bfe(\xi)=0$ so the
constant densities are excluded.  The maximum is achieved
at multiplies and rotations of the field $\ell(\xi)$ 
from the next definition.

\begin{definition}
\label{def:ell}  For $\xi\in S^{d-1}$ denote by $\ell(\xi)$
the projection of the vector $(1,0,\dots,0)$
orthogonal to $\xi$,
\begin{equation}
\label{eq:defell}
\ell(\xi)  
\ :=\
(1,0,\dots,0)
\ -\ 
\big(\xi\cdot(1,0,\dots,0)\big)\,\xi 
\ =\ 
(1,0,\dots,0)
\ -\ 
\xi_1\,\xi\,.
\end{equation}
\end{definition}

$\ell$ is a vector field whose integral curves are the lines
of longitude connecting the pole $(-1,0,\cdots,0)$ to
opposite pole
$(1,0,\cdots ,0)$.
The maximum value of $J_1$ for electromagnetic
waves is smaller by $(d-1)/d$ 
than the extremum in the scalar case.   
The same functions also
solve the $J_2$ problem
when $R$ is not too large.
The study of $J_1$ is reduced to an application of 
the Cauchy-Schwartz inequality.   

In \S \ref{sec:equiv}
the maximization
of $J_2$ is transformed to a 
problem in spectral theory.
Maximizing $J_2$ 
 is equivalent to finding the norm of
an operator.  In the scalar case we call the operator
$L$.   Finding the norm  is equivalent to finding
the spectral radius of the self adjoint operator $L^*L$.
The operator $L^*L$ is compact and rotation invariant
on $L^2(S^{d-1})$.  Its spectral theory is reduced
by the spaces of spherical harmonics of order $k$.  
On the space of spherical harmonics
of degree $k$,
$L^*L$
 is multiplication 
by a constant 
$
\Lambda_{d,k}(R)$ 
computed exactly in terms of Bessel functions
in 
Theorem \ref{thm:lambdaofR}.
Theorem 7.1 shows that 
in the scalar case $\Lambda_{d,0}(R)$ is 
the largest for $R\le \pi/2$.

For the Maxwell problem, the corresponding operator
is 
denoted $\bfL_M^*\bfL_M$.   We do not know all its
eigenvalues.   However two explicit eigenvalues
are $(2/3)(\Lambda_{d,0}(R)-\Lambda_{d,2}(R))$,
and, $\Lambda_{d,1}(R)$.
When $R\le \pi/2$,
Theorem 7.3 proves that  the they are the largest and
second largest
eigenvalues.  The proof uses the minimax principal.
 The eigenfunctions 
for the largest eigenvalue are rotates of multiplies of $\ell$.

 For $R>\pi/2$
we  derive 
in \S 8
rigorous sufficient conditions guaranteeing
that the 
same functions provide the extremizers.
The conditions involve the $\Lambda_{d,k}$.
To verify them we evaluate the integrals defining the
$\Lambda_{3,k}$ approximately.
By such evaluations we show  that 
for $d=3$ and $R\le 2.5$ the solutions maximizing
$J_1$ 
also maximize $J_2$.

The energy density  is equal to  the largest eigenvalue
divided by $|B_R(0)|$.
In the range $d=3, R\le 2.5$ 
in 
both the scalar and electromagnetic case 
this quantity is a decreasing
functions of $R$.  This shows that 
the third problem at the start, of finding the radius
with highest energy density is solved by $R=0$.
However the graph is fairly flat.
The density dips to about 1/2 its maximum
at about $R=2$ which is about a third the wavelength.

For 
focusing of electromagnetic waves to a ball of 
radius $R$ no larger than one third of a wavelength the optimal
strategy is to choose $\bfe(\xi)$ a multiple of a rotate of 
$\ell(\xi)$.  The extremizing electric fields
are as polarized as a divergence free field can be. 
When $d=3$ formula \eqref{eq:stationary} and Example
\ref{ex:incoming} 
show that the far field for this choice is equal up to rotations by
$$
c\
\frac{\sin|x|}{|x|}
\
\ell\Big(
\frac{ x }{  | x | }
\Big)\,.
$$
This field is cylindrically symmetric with axis of symmetry
along the $x_1$-axis.   
The restriction of $\ell$ to the unit sphere is cylindrically symmetric
given by rotating the following figure about the horizontal axis.

\begin{figure}
\begin{center}
\includegraphics[width=0.5\textwidth]{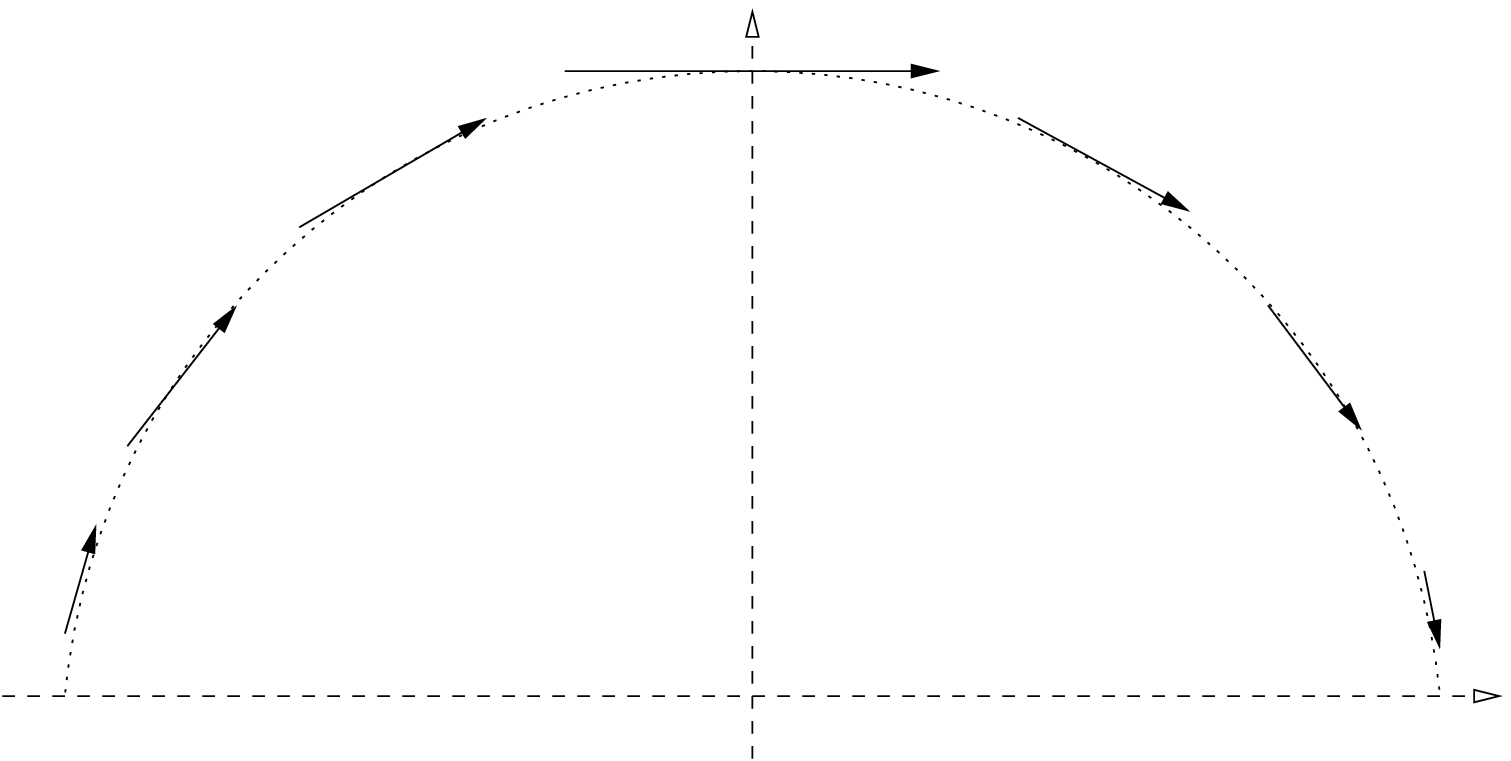}
\end{center}
\end{figure}

For the problem of focusing a family of lasers, 
this suggests using linearly polarized 
sources concentrated near vertical equator and sparse near the
poles on the horizontal axis. In contrast, for scalar waves one should distribute
sources as uniformly as possible.  

\vskip.2cm

{\bf Acknowlegements.}   I  thank G. Mourou for 
proposing this problem.  
I early conjectured that the constants and $\ell(\xi)$ were
the extremizers in the scalar and electromagnetic cases
respectively.
J. Szeftel, G. Allaire,
 C. Sogge, P. G\'erard, and J. Schotland
provided both encouragement and help on the path
to the results presented here.
The meetings were in  Paris, Pisa, and Lansing. 
In Europe I was a guest at  the Ecole
Normale Sup\'erieure, Universit\'e de Paris Nord, and,
Universit\`a di Pisa.
Sogge
and G\'erard were guests of  the Centro De Giorgi.
Schotland and I were both guests
of the IMA and Michigan State University.
I thank all these individuals and institutions.

\section{Monochromatic waves.}

\subsection{Electromagnetic waves and their transforms.}

\begin{proposition}
\label{prop:divergencefree}
{\bf i.}  $E$ given by \eqref{eq:E} satisfies
$
\dive E =0
$
if and only if
\begin{equation}
\label{eq:Fperp}
\bfe(\xi)\cdot \xi \ =\ 0,
\qquad
{\rm on}
\qquad
\{|\xi|=1\}\,.
\end{equation}

{\bf ii.}
If a monochromatic solution of the 
Maxwell equations has electric field
given by \eqref{eq:monoplusminus} with $\omega=1$
 and $E$ is given by \eqref{eq:E} then the magnetic
 field  is equal to $e^{ it}\,B(x)$ with 
\begin{equation}
\label{eq:B}
B(x)
\ =\
-\
\int_{|\xi|=1}
e^{ix\xi}\
\xi\wedge \bfe(\xi)\
d\sigma\,,
\end{equation}
\end{proposition}

{\bf Proof.}   
Differentiating \eqref{eq:E} yields,
$$
\dive E
\ =\ 
\int_{|\xi|=1}
e^{ix\xi}\
i\,\xi\cdot \bfe(\xi) d\sigma
\,,
\qquad
\curl E \ =\ 
\int_{|\xi|=1}
e^{ix\xi}\
i\,\xi\wedge \bfe(\xi)\
d\sigma\,.
$$
The first formula proves {\bf i.}    

The Maxwell equations together with 
\eqref{eq:monoplusminus}
yield
$$
-\,\curl E \ =\ B_t \ =\ 
 \,i\, B\,.
$$
Therefore, the second formula proves {\bf ii.}
\qed

\vskip.2cm

\begin{remark}
The condition \eqref{eq:Fperp} asserts that 
$\bfe(\xi)$ is tangent to the unit sphere.   
Brouwer's Theorem asserts that  if $\xi\mapsto\bfe(\xi)$ is continuous
then there must be a $\underline\xi$ where $\bfe(\underline\xi)=0$.
\end{remark}

\begin{example}.  If $d=3$ and $E$ is given \eqref{eq:E} 
with $\bfe(\xi)=\ell(\xi)$
then on $|\xi|=1$,
$$
\xi\wedge \ell(\xi)
\  =\
\xi\wedge\Big((1,0,0)\ -\ \xi_1\,\xi\Big)
\ =\
\xi\wedge(1,0,0)
\ =\
(0, \xi_3, -\xi_2)
$$
is the tangent field to latitude lines winding
around the $x_1$-axis.   Since this is an odd
function,     the magnetic field vanishes at the origin, $B(0)=0$.
 \end{example}

\subsection{The competing solutions.}
\label{sec:competing}

First verify the stationary phase formula \eqref{eq:stationary}
from the appendix.
Consider  $E$ given by \eqref{eq:E} with $\bfe\in C^\infty(\{|\xi|=1\})$.
For $x$ large, the integral  \eqref{eq:E} has two
stationary points, $\xi=\pm x/|x|$.  At $\xi=x/|x|$ 
parameterize the surface by coordinates in the tangent
plane at $x$ to find that 
the phase
$x \xi$ has a strict maximum equal to $|x|$ and hessian
equal to the $-I_{(d-1)\times(d-1)}$.  At $\xi =
-x/|x|$ the phase has a minimum with value $-|x|$ and  hessian
equal to the identity.
The stationary phase method yields
\eqref{eq:stationary}.
The energy in the electric field
satisfies \eqref{eq:slowgrow} and \eqref{eq:L2}.

\begin{theorem}
Suppose that $E\in {\cal S}^\prime(\RR^d)$
is a tempered solution of the reduced wave
equation given by \eqref{eq:E} with $\bfe\in 
{\cal D}^\prime(S^{d-1})$.  Then,
\eqref{eq:slowgrow} holds 
if and only if 
$\bfe \in L^2(S^{d-1};\CC^d)$.  In that case
\eqref{eq:L2} holds.
In addition the stationary phase
approximation holds in the sense that
as $R\to\infty$,
$$
\int_{R\le |x|\le 2R}
\Big| 
E(x)  -
\frac{1/\sqrt{2\pi}}
{
|x|^{(d-1)/2}
 }
\
\Big(
e^{-i | x |}\,\bfe((-x/|x|)
+
e^{i\pi(d-1)/4}\, e^{i|x|}\,
\bfe(x/|x|)
\Big)
\Big|^2\,dx
= o(R)\,.
$$
\end{theorem}

{\bf Proof.}   The first two assertions are consequences of 
H\"ormander \cite{Hor} Theorems 7.1.27
and 7.1.28.  
The Theorem 7.1.28 also implies that
$$
\sup_{R\ge 1}\ \frac{1}{R}
\int_{R\le |x|\le 2R}
\big| 
E(x) 
\big|^2\ dx
\ \le \
c(d) \,  \int_{|\xi | = 1 }
\big|
\bfe(\xi)
\big|^2\ d\sigma\,.
$$
This estimate shows that to prove the third assertion
it suffices to prove it for the dense set of $\bfe\in C^\infty(S^{d-1})$.
In that case the result is a consequence of the stationary phase
formula \eqref{eq:stationary}.
\qed

\vskip.2cm

\begin{definition}
$\bfH$ is the closed subspace  of $\bfe\in L^2(S^{d-1};\CC^d)$
consisting of $\bfe(\xi)$ so that $\xi\cdot\bfe(\xi)=0$. 
Denote by $\Pi$ the orthogonal projection 
of 
$\bfe\in L^2(S^{d-1};\CC^d)$
on $\bfH$.
\end{definition}

The next example explains a connection
between the solutions of the reduced equation that
we consider and those satisfying the Sommerfeld
radiation conditions.

\begin{example}  If $g\in {\cal E}^\prime(\RR^d:\CC^d)$
is a distribution with compact support, then there
are unique solutions of the reduced wave equation
$$
(\Delta + 1)E_{out}\ =\ g,
\qquad
({\rm resp.}\quad
(\Delta + 1)E_{in}\ =\ g)
$$
satsifying the outgoing (resp. incoming) radiation
conditions.
The difference $E:=E_{out}-E_{in}$ is a solution
of the homogeneous reduced wave equation.
The field $F:=e^{it}E(x)$ is the unique solution
of the initial value problem
$$
\Box F\ =\ 0,
\qquad
F\big|_{t=0}\ =\ g,
\quad
F_t\big|_{t=0}\ =\ ig\,.
$$
When $g\in L^2$ 
 the formula $E(x)\delta(\tau-1)=c\int_{-\infty}^\infty e^{-i\tau t} F\,dt$
 together with the solution formula for the Cauchy problem
imply  that
\eqref{eq:slowgrow}
holds
 (or see \cite{Hor} Theorem 14.3.4 showing that both incoming
 and outgoing fiels satisfy \eqref{eq:slowgrow}).  More generally,
 the Fourier transforms in time of solutions of Maxwell's
 equations with compactly supported divergence free
 square integrable
  initial data yield examples of  monochromatic solutions in our class.
   \end{example}

\subsection{Spherical symmetry is impossible.}

It is natural to think that focusing is maximized
if waves come in equally in all directions.
For the scalar wave equation that is the case.
  However,
such waves do not exist
for Maxwell's equations.
Whatever is the definition of spherical symmetry,
such a  field must satisfy the hypotheses
of the following theorem.

\begin{theorem}  If $E(x)\in C^1(\RR^d)$ satisifies 
$\dive E=0$ and for $x\ne 0$ the angular part
$$
E
\ -\ 
\bigg(E\cdot \frac{x}{ | x | }\bigg)\,\frac{x}{ | x | }
$$
has length that depends only on $|x|$ then
$E$ is identically equal to zero.
\end{theorem}

{\bf Proof.}  The restriction of  the angular part of $E$ to each sphere
$|x|=r$ is a $C^1$ vector field tangent to the sphere and of constant length.
Brouwer's Theorem asserts that there is a point $\ux$ on the sphere
where the tangent vector field vanishes.  Therefore the constant length
is equal to zero and $E$ is radial.

Therefore in $x\ne 0$,
$$
E(x) = \phi(|x|)\, x
\,.
$$
Since $E\in C^1$ it follows that $\phi\in C^1(\{|x|>0\})$.

Compute for those $x$,
$$
\dive 
E\ =\ \phi \, \dive x \ +\ 
(\nabla_x\phi) \cdot x
\ =\ 
d\,\phi  \ +\  r\,\phi_r \,.
$$
Therefore in $x\ne 0$, $r\,\phi_r+d\,\phi=0$ so $\phi=c\, r^{-d}$.

Since $E$ is continuous at the origin it follows that $c=0$ so 
$E=0$ in $x\ne 0$.  By continuity, $E$ vanishes identically.
\qed

\section{Maximum field strengths.}

We
solve the variational problems associated to
the functional $J_1$ to yield sharp
pointwise bounds on monochromatic waves.
The fact that the bounds for electromagnetic
fields are smaller shows that
focusing effects are  weaker.
The extremizing fields are first
characterized by their Fourier Transforms.
Explicit formulas in $x$-space 
are given in
\S \ref{subsec:4.4}

\subsection{Scalar waves.}

\begin{theorem}
\label{thm:scalarmaxpoint}
If 
$$
u(x) \ = \ 
\int_{|\xi|=1}
e^{ix\xi}\ f(\xi)\ d\sigma\,,
\qquad
f
\ \in\ L^2(S^{d-1} ),
$$
and $\ux\in \RR^d$,
then
$$
\big| u(\ux) \big|
\ \le \
|S^{d-1}|^{1/2}
\,
\big\| f  \big\|_{L^2(S^{d-1})}
$$
with equality achieved if an only if $f$ is a scalar multiple of $e^{-i\ux\xi}$.
\end{theorem}

{\bf Proof.}  The quantity to maximize is the $L^2(S^{d-1})$ scalar
product of $f$ with $e^{-ix\xi}$.  The result is exactly the 
Cauchy-Schwartz inequality.
\qed

\subsection{Electromagnetic waves.}

From Definition \ref{def:ell}, 
 $\ell(\xi)$ is tangent to the longitude 
lines on the unit sphere connecting the pole $(-1,0,\dots,0)$ to the
pole $(1,0,\dots,0)$.  It is the gradient of the restriction of the 
function $\xi_1$ to the unit sphere.

\begin{theorem}
\label{thm:maxpoint}
If $d\ge 2$ and
$$
E(x)=
\int_{|\xi|=1} e^{ix\xi}\ \bfe(\xi)\ d\sigma\,,
\qquad
\bfe\ \in\
\bfH\,.
$$
Then
\begin{equation}
\label{eq:Ebound}
| E(0) | \ \le\ 
\bigg(
\frac{d-1}{d}\ 
|S^{d-1}|
\bigg)^{1/2}\
\big\| \bfe \big\|_{L^2(S^{d-1})}\,.
\end{equation}
Equality holds if and only if $\bfe$ is equal to a constant
mulitple of a rotate of 
$\ell(\xi)$.
\end{theorem}

{\bf Proof.}  By homogeneity it suffices to consider
$\|\bfe\|_{L^2(S^{d-1})}=1$.
Rotation and multiplication by a complex
number of modulus one
reduces to the case $E(0)=|E(0)|(1,0,\dots,0)$
and  $ |E(0) | = \int \bfe_1(\xi)\,d\sigma$.
Need to study,
 $$
\sup\ \Big\{  \int \bfe_1(\xi)\, d\sigma\ :\ 
\ \xi.\bfe(\xi)=0, \ \ \int_{|\xi|=1} \|\bfe(\xi)\|^2\,d\sigma \ =\ 1\Big\}\,.
$$
The quantity to be maximized is
$$
\int_{| \xi|=1}
\bfe_1(\xi)\,d\sigma
\ =\
\big(
\bfe\,,\,
(1,0.\cdots,0)
\big)_{L^2(S^{d-1})}
$$
The constant function,
$(1,0.\cdots,0)$ does not belong to the subspace
$\bfH$.  The projection theorem shows that the quantity is
maximized for $\bfe$ proportional to the projection of 
$(1,0.\cdots,0)$ on $\bfH$.  Equivalently,
using \eqref{eq:defell} together with $\bfe\cdot\xi=0$ yields
$$
\bfe_1
\ =\ 
\bfe\cdot (1,0,\dots,0)
\ =\ 
\bfe\cdot
\big(
\ell(\xi) \ +\ \xi_1\,\xi\big)
\ =\ 
\bfe\cdot \ell\,.
$$

\begin{equation}
\label{eq:inte1}
\int_{| \xi|=1}
\bfe_1(\xi)\,d\sigma
\ =\ 
\int_{| \xi|=1}
\bfe(\xi)\cdot \ell(\xi)\ d\sigma
\end{equation}

which is equal to the $\bfH$ scalar product of 
$\bfe$ and $\ell$.  Since one has the orthogonal
decomposition
$$
(1,0,\dots , 0)
\ =\ 
\ell(\xi) + \xi_1\,\xi,
\qquad
{\rm one\ has},
\qquad
1\ =\ 
|\ell(\xi)|^2  \ +\ \xi_1^2\,.
$$
The Cauchy-Schwartz inequality
shows that the quantity
\eqref{eq:inte1} is
\begin{equation}
\label{eq:CS}
 \le \
\|\bfe \|_{L^2(S^{d-1})}
\,
\|\ell \|_{L^2(S^{d-1})}
\  =\   
\|\bfe \|_{L^2(S^{d-1})}
\,
\bigg(\int_{|\xi|=1} (1-\xi_1^2)\ d\sigma\bigg)^{1/2}
.
\end{equation}
The extremum is attained uniquely
when  $\bfe =z\,\ell/\|\ell\|$ with $|z|=1$.  

To evaluate the 
integral on the right of \eqref{eq:CS} compute,
$$
\int_{|\xi|=1} \xi_1^2\, d\sigma
\ =\
\int_{|\xi|=1} \xi_j^2\, d\sigma
\ =\ 
\frac{1}{d} \int_{|\xi|=1} \sum_j\xi_j^2\ d\sigma
\ =\ 
\frac{1}{d}\int_{|\xi|=1} 1\,d\sigma
\ =\
 \frac{|S^{d-1}|}{d}\,.
$$
Therefore
$$
\int_{|\xi|=1} (1-\xi_1^2)\ d\sigma
\ =\ 
|S^{d-1}| \ -\
 \frac{|S^{d-1}|}{d}
 \ =\ 
 |S^{d-1}| \
 \frac{ d-1 } { d } \,.
 $$
 Together with \eqref{eq:CS} this proves
\eqref{eq:Ebound}.
\qed

\begin{remark}   If one constrains $\bf e$ to have
support in a subset $\Omega$ then with $\chi$ denoting
the characteristic function of $\Omega$,
\begin{equation*}
\int_{|\xi|=1}
\bfe(\xi)\cdot\ell(\xi)\ d\sigma
=
\int_{|\xi|=1}
\bfe(\xi)\cdot\ell(\xi)\chi(\xi)\ d\sigma
\end{equation*}
and $E_1$ is maximized by the choice
$\bfe = \ell(\xi)\chi(\xi)$.   In the extreme light initiative
$\Omega$ is a small number of disks distributed
around the equator  $x_1=0$.  
\end{remark}

\subsection{Derivative bounds.}
\label{sec:derivatives}

\begin{corollary} If $d\ge 2$ and $E$ satisfies 
$$
E(x)=
\int_{|\xi|=1} e^{ix\xi}\ \bfe(\xi)\ d\sigma\,,
\qquad
\bfe\ \in \
\bfH\,,
$$
then for all
$\alpha\in \NN^d$ and $\ux\in \RR^d$,
\begin{equation}
\label{eq:DEbound}
\big|
\partial_x^\alpha E(\ux) 
\big|
\ \le\
\bigg(
\frac{d-1}{d}\ 
|S^{d-1}|
\bigg)^{1/2}
\
\big\|\bfe\big\|_{L^2(S^{d-1})}
\,.
\end{equation}
\end{corollary}

{\bf Proof.}  The case $\alpha=0$ follows from Theorem
\ref{thm:maxpoint} applied to 
$$
\widetilde E(x) 
\ :=\
 E(x+\ux)
\ = \
\int_{|\xi|=1}
e^{ix\xi}
\ 
e^{i\ux\xi}\,
\bfe(\xi)\
d\sigma
\ :=\
\int_{|\xi|=1}
e^{ix\xi}
\ 
\widetilde \bfe(\xi)\
d\sigma
\,.
$$

Compute for $|\alpha|>0$
$$
\partial_x^\alpha E \ = \
\partial_x^\alpha
\int_{|\xi|=1}
e^{ix\xi}\, \bfe(\xi)\ d\sigma
\ =\
\int_{|\xi|=1}
e^{ix\xi}\
(i\xi)^\alpha\bfe(\xi)\ d\sigma
\,,
$$
which is of the same form as $E$ with 
density $(i\xi)^\alpha\bfe(\xi)$ orthogonal
to $\xi$.  Since  $|\xi_j|\le 1$ 
it follows that $|\xi^\alpha|\le 1$ so  
$
\|(i\xi)^\alpha\bfe\|_{L^2(S^{d-1})}
\le \| \bfe\|_{L^2(S^{d-1})}
$.
Therefore the general case 
follows from the case $\alpha=0$.
\qed

\begin{remark}    Derivative bounds for the scalar case
are derived in the same way.  They lack the factor 
$(d-1)/d$.
\end{remark}

\subsection{Formulas for the extremizing fields.}
\label{subsec:4.4}

The electric 
field corresponding to the extremizing
density $\ell$ is explicitly calculated.  
The computation  relies on 
relations between Bessel functions, spherical harmonics,
and, the Fourier Transform.  These relations
are needed to analyse $J_2$.

Start from 
identities in   Stein-Weiss \cite{Stein}.
Their Fourier transform is defined on 
page 2,
$$
\int f(x)\ e^{-i2\pi x\xi}\ dx,
\qquad
{\it n.b.}   {\rm \ the \ } 2\pi{\rm\ in\ the \ exponent.}
$$ 
We will not follow this convention, so  
adapt
their identities.
The Bessel function of order $k$ is 
(page 153),
\begin{equation}
\label{eq:Jdef}
J_k(t) =
\frac{(t/2)^k}{\Gamma[(2k+1)/2]\ \Gamma(1/2)}\
\int_{-1}^1
e^{its}\
(1-s^2)^{(2k-1)/2}\
ds,
\ \ 
-1/2 <k\in  \RR\,.
\end{equation}

Theorem 3.10 (page 158) is the following.

\begin{theorem}
If $x\in \RR^d$, $f= f(|x|)P(x)\in L^1(\RR^d)$ with $P$ a homogeneous
harmonic polynomial of degree $k$, 
then $\int f(x)e^{-2\pi i x \xi}dx=F(|\xi|)P(\xi)$
with
$$
F(r) \ =\ 
2\pi\,
i^{-k}\,
r^{-(d+2k-2)/2}\,
\int_0^\infty
f(s)\
J_{(d+2k-2)/2}(2\pi r s)\
s^{(d+2k)/2}\
ds\,.
$$
\end{theorem}

This theorem is equivalent, by scaling and linear combination,
to the same formula with  $f=\delta(r-1)$.
That case is the  identity,
\begin{equation}
\label{eq:steinP}
\int_{|x|=1}
e^{-i2\pi x\xi}\
P(x)\
d\sigma
\ =\
2\pi\,
i^{-k}\,
|\xi|^{-(d+2k-2)/2}\,
J_{(d+2k-2)/2}(2\pi  |\xi|)\
P(\xi).
\end{equation}

\begin{remark}  {\bf i.} For $|\xi|\to \infty$, $J(|\xi|)=O(|\xi|^{-1/2})$,
and 
$P(\xi)=O(|\xi|^k|)$ so the right hand side is 
$O(|\xi|^{-(d-2)/2-1/2})= O(|\xi|^{-(d-1)/2})$
as required by the principle of stationary phase.

{\bf ii.}
For $|\xi|\to 0$, $J_{((d-2k-2)/2}(|\xi|) = O(|\xi|^{(d-2k-2)/2})$
so the right hand side of \eqref{eq:steinP} is $O(|\xi|^k)$.
The higher the order of $P$ the smaller is the Fourier
transform near the origin.
\end{remark}

To adapt to the Fourier transform without the $2\pi$ in 
the exponent,
use the substitution $\eta=2\pi \xi$, $|\eta|=2\pi |\xi|$
 to find,
  \begin{equation*}
\label{eq:hormP}
\int_{|x|=1}
e^{-ix\eta}\
P(x)\
d\sigma
\ =\
2\pi\,
i^{-k}\,
(|\eta|/2\pi)^{-(d+2k-2)/2}\,
J_{(d+2k-2)/2}(|\eta|)\
P(\eta/2\pi)\,.
\end{equation*}
Using the homogeniety of $P$ yields
$$
\ =\
(2\pi)^{1-k}\,
i^{-k}\,
(|\eta|/2\pi)^{-(d+2k-2)/2}\,
J_{(d+2k-2)/2}(|\eta|)\
P(\eta)\,.
$$
The exponent of $2\pi$ is equal to
$d/2$
yielding,
\begin{equation}
\label{eq:harmonic2}
\int_{|x|=1}
e^{-ix\eta}\
P(x)\
d\sigma
=
(2\pi)^{d/2}\,
i^{-k}\,
|\eta|^{-(d+2k-2)/2}\,
J_{(d+2k-2)/2}(|\eta|)\
P(\eta).
\end{equation}
Since, $|\eta|^{-k}P(\eta)=P(\eta/|\eta|)$ 
\eqref{eq:harmonic2} equivalent to,
\begin{equation}
\label{eq:FThomog} 
\int_{|x|=1}
e^{-ix\eta}\
P(x)\
d\sigma
=
(2\pi)^{d/2}\,
i^{-k}\,
|\eta|^{-(d-2)/2}\,
J_{(d+2k-2)/2}(|\eta|)\
P(\eta/|\eta|).
\end{equation}
The change of variable $\eta\mapsto -\eta$ yields,
\begin{equation}
\label{eq:FThomog2} 
\int_{|x|=1}
e^{ix\eta}\
P(x)\
d\sigma
=
(2\pi)^{d/2}\,
(-i)^{-k}\,
|\eta|^{-(d-2)/2}\,
J_{(d+2k-2)/2}(|\eta|)\
P(\eta/|\eta|).
\end{equation}
Finally interchange the role of $x$ and $\eta$ to find,
\begin{equation}
\label{eq:FThomog3} 
\int_{|\eta|=1}
e^{ix\eta}\
P(\eta)\
d\sigma
=
(2\pi)^{d/2}\,
(-i)^{-k}\,
|x|^{-(d-2)/2}\,
J_{(d+2k-2)/2}(|x|)\
P( x / |x| ).
\end{equation}

\begin{example}
\label{ex:scalmax}
  The second most interesting example
is the extremizing field for the scalar case when $d=3$.  In that case 
$P=constant$ and there
is a short derivation.  The function $u(x):=\int_{| \xi | =1 } e^{ix\xi}\,d\sigma$
is a radial solution of $(\Delta+1)u=0$.  In $x\ne 0$ these
are spanned for $d=3$ by $e^{\pm ir}/r$.  Smoothness
at the origin forces $u=A\sin r/r$.  Since
$u(0)=|S^{d-1}|$
it follows that
$A=|S^{d-1}|$.
\end{example}

 The most interesting case for us is $d=3$ and
  the extremizing field $E$ 
with $\bfe(\xi)=\ell(\xi)$.
Since $\ell$ is not a spherical harmonic, the preceding
result does not apply directly.   
To find the exact electric field, decompose $\ell$
in spherical harmonics.  

 \begin{lemma}
  \label{lem:ellexp}
  The
   spherical harmonic expansion of the restriction of $\ell(\xi)$
   to the  unit sphere $S^{d-1}\subset \RR^d$ is 
   \begin{equation}
   \label{eq:ellexp}
 \ell(\xi)\ =\  \bigg(
   \frac{d-1}{d} -
   \sum_{j=2}^d
   \frac{\xi_1^2-\xi_j^2}{d}
\,,\,
-\xi_1\xi_2
\,,\,
\cdots
\,,\,
-\xi_1\xi_d
\bigg)
\quad
{\rm on}
   \quad
   |\xi|=1\,.
\end{equation}
\end{lemma}

  {\bf Proof of lemma.}
  When $|\xi |=1$, 
   \begin{equation}
   \label{eq:ellstep1}
   \ell(\xi)
   =
   (1,0,\dots,0) -
   \xi_1\xi
    = 
   (1,-\xi_1\xi_2, \dots, -\xi_1\xi_d) -  (\xi_1^2,0,\dots,0)
   \,.
   \end{equation}
   The first summand has coordinates that are
    spherical harmonics. 
  
  Decompose
  $$
 \xi_1^2
   \ =\ 
   \frac{\xi_1^2 + \dots + \xi_d^2}{d}
   \ +\
   \sum_{j=2}^d
   \frac{\xi_1^2-\xi_j^2}{d}\,,
   \qquad
   \xi\in \RR^d\,,
   $$ 
   to find the expansion in spherical harmonics of the restriction of $\xi_1^2$
   to the unit sphere
   $\xi_1^2 + \dots + \xi_d^2=1$,
   \begin{equation}
   \label{eq:xi1square}
   \xi_1^2
    \ =\ 
     \frac{1}{d}
   \ +\
   \sum_{j=2}^d
   \frac{\xi_1^2-\xi_j^2}{d}
   \qquad
   {\rm on}
   \qquad
   \xi_1^2 + \cdots +  \xi_d^2=1\,.
   \end{equation}

   Using \eqref{eq:xi1square} in \eqref{eq:ellstep1} proves
   \eqref{eq:ellexp}.
   \qed

\begin{example}
\label{ex:incoming}  The stationary phase formula
\eqref{eq:stationary} applied to the extremizing
$\bfe=\ell(\xi)$ which is an even function yields
for $d=3$
$$
E(x)
=
\frac{-1}{\sqrt{2\pi} }\
\ell\Big(
\frac
{x}
{ |x| }\Big)\
\frac
{\sin |x|}
{ |x| }
\ +\ O(|x|^{-2})
\,.
$$
This is the sum of  incoming
and outgoing waves with spherical wave fronts and each with
profile on large spheres proportional to $\ell(x/|x|)$.  The desired
incoming wave is such an $\ell$-wave.  
\end{example}

\section{Equivalent selfadjoint eigenvalue problems.}
\label{sec:equiv}

The section introduces eigenvalue problems equivalent
to the maximization of $J_2$.

\subsection{The eigenvalue problem for focusing scalar waves.}

\begin{definition}  For $R>0$ 
define the  compact linear operator $L:L^2(S^{d-1})\to L^2(B_R(0))$
by
$$
(Lf)(x)
\ :=\
\int_{ |\xi|=1} e^{ix\xi}\ f(\xi)\ d\sigma
\,.
$$
\end{definition}

The operator $L$ commutes with rotations. The adjoint
$L^*$ maps $L^2(B_R(0))\to L^2(S^{d-1})$.  

\begin{proposition}  The following four problems are equivalent.

{\bf i.}  Maximize the functional $J_2$ on scalar monochromatic
waves.

{\bf ii.}  Find $f\in L^2(S^{d-1})$ with $\|f\|_{L^2(S^{d-1} ) } =1$
so that $\|Lf\|_{B_R(0)}$ is largest.

{\bf iii.}  Find
the norm of $L$.

{\bf iv.}  Find the largest eigenvalue
of the positive compact self adjoint operator $L^*L$
on $L^2(S^{d-1})$.
\end{proposition}

{\bf Proof.}   The equivalence of the first three follows from the 
definitions.
The equivalence with the third follows from the identity
$$
\|Lf\|^2_{L^2(B_R)}
\ =\ 
(Lf\,,\, Lf)_{L^2(B_R)}
\ =\ 
(L^*Lf\,,\, f)_{L^2(S^{d-1})}
\,.
\eqno{
\qed}
$$

\vskip.2cm

\begin{definition}
Define a rank one operator
$$
 L^2(S^{d-1};\CC)
\ \ni\ 
f
\  \to \ 
L_0f\ :=\ 
\int_{| \xi | = 1}
f(\xi)\ d\sigma
\ \in \ \CC
\,.
$$
\end{definition}

\begin{remark}  {\bf i.} The problem of maximizing $J_1$ for scalar waves is
equivalent to finding the norm of $L_0$ and also finding
the largest eigenvalue of $L_0^*L_0$.

{\bf ii.}  The same formula defines an operator from $L^2(S^{d-1})
\to L^2(B_R(0))$ mapping $f$ to a constant function.  With only small
risk of confusion we  
use the same symbol $L_0$ for that operator too.
\end{remark}

\begin{definition}  The vector valued version of  $L$ and $L_0$ are 
defined by
$$
(\bfL  \bfe)(x)
\ :=\ 
\int_{| \xi | = 1}
e^{ix\xi}\ \bfe(\xi)\ d\sigma\,,
\qquad
\bfe\in L^2(S^{d-1};\CC^d)\,,
$$
$$
\bfL_0\bfe\ :=\ 
\int_{| \xi | = 1}
\bfe(\xi)\ d\sigma\,,
\qquad
\bfe\in L^2(S^{d-1};\CC^d)\,.
$$
 Denote by $\bfL_M$ and 
$\bfL_{0,M}$ the 
restriction 
of $\bfL$ and $\bfL_0$ to $\bfH$.
\end{definition}

\begin{remark}
\label{rmk:LL}
 {\bf i.}  The problem of maximizing the 
functional $J_1$ for monochromatic electromagnetic
waves is equivalent to finding the largest eigenvalue
of $\bfL_{0,M}^*\bfL_{0,M}$.

{\bf ii.}  The problem of maximizing
the functional $J_2$ for 
monochromatic electromagnetic
waves is equivalent to finding the 
largest eigenvalue of 
$\bfL_M^*\bfL_M$.

{\bf iii.}  The operator $\bfL \Pi$ is equal to $\bfL_M$ on $\bfH$
and equal to zero on $\bfH^\perp$.  Therefore
$(\bfL \Pi)^*(\bfL\Pi)$ is equal to $\bfL_M^*\bfL^{}_M$ on $\bfH$
and equal to zero on $\bfH^\perp$.
\end{remark}

\section{Exact eigenvalue computations.}

\subsection{The operators $L_0^*L_0$, $\bfL_0^*\bfL_0$ and $\bfL_{0,M}^*\bfL_{0,M}$.}

\begin{theorem}
\label{thm:Lzerospec}   {\bf i.} The spectrum of  $L_0^*L_0$ contains 
one nonzero eigenvalue, $|S^{d-1}|$, with multiplicity one.
The eigenvectors are the constant functions.

{\bf ii.}  The spectrum of the operator $\bfL_0^*\bfL_0$ contains
one nonzero eigenvalue,
$|S^{d-1}|$, with multiplicity $d$.
The eigenvectors are the $\CC^d$-valued constant functions.

{\bf iii.}  The spectrum of $\bfL_{0,M}^*\bfL_{0,M}$ 
contains one nonzero
eigenvalue, $
|S^{d-1}|(d-1) 
/d
$
with multiplicity $d$.  The corresponding eigenspace
consists of 
$\ell(\xi)$ and functions obtained by rotation and scalar
multiplication.
\end{theorem}

 {\bf Proof iii.}  Suppose that $f$ is an eigenfunction in $\bfH$
 so that the norm of $\int f\ d\sigma$ is maximal.
 Rotating $f$ yields a function in $\bfH$ with the same
 $\|\bfL_0 f\|$ and with 
 $\bfL_0 f$ parallel to $(1,0,\dots,0)$.  Therefore maximizing
 $\bfL_0 f$ and maximizing
 $$
 \bfL_0 f \cdot (1,0,\dots\,,0)
 $$
 yield the same extreme value.  
 
 Since $f\in \bfH$,
 $$
 \int f_1 d\sigma
 \ =\ 
 \int f\cdot (1,0,\dots,0))\ d\sigma
 \ =\ 
 \int f\cdot \ell\ d\sigma\,.
 $$
 The extreme value is attained for $f$ parallel to $\ell$.
 This shows that $\ell$ is an eigenfunction corresponding
 to the largest eigenvalue.

 Rotating and taking scalar multiples yields a complex 
 eigenspace of dimension $d$.
  Since ${\rm rank}\,\bfL_0=d$ these are all the eigenfunctions.
 
If $\|f\|_{L^2(S^{d-1} ) } =1$, then the maximization
 of $J_0$ shows that  
$|\bfL_{0,M}f|^2 =
|S^{d-1}|(d-1) 
/d
$ 
proving the formula for the eigenvalue.
  \qed
  
  \subsection{The eigenfunctions and eigenvalues of $L^*L$ and $\bfL^*\bfL$.}

\begin{theorem}
\label{thm:lambdaofR}
In dimension $d$, the spherical harmonics of order
$k$ 
are eigenfunctions of $L^*L$ 
with eigenvalue
\begin{equation}
\label{eq:Lambdadef}
\Lambda_{d,k}(R)\ :=\
(2\pi)^{d/2}\,
|S^{d-1}|
\ 
\int_0^R
r
\ 
\big[
J_{(d+2k-2)/2}(r)
\big]^2
\ dr\,.
\end{equation}
\end{theorem}

\begin{remark}
  From the point of view of focusing of 
energy 
into balls, all spherical harmonics of the same order
are equivalent.  
\end{remark}

{\bf Proof.}  If $P$ is a homogeneous harmonic
polynomial of degree $k$ formula  
\eqref{eq:FThomog3} 
shows that 
$$
(L\,P)(x) \ =\ 
\phi_{d,k}( | x | )\
P(x),
$$
defining the function $\phi$.  

The operator $L^*$ is an integral operator from
$L^2(B_R)\to L^2(S^{d-1})$ with kernel
$e^{-ix\xi}$.  Therefore
$$
L^*L(P)
\ =\ 
\int_{|x|\le R}
e^{-ix\xi}
\
\phi_{d,k}( | x | )\
P(x)
\
dx\,.
$$
Introduce polar coordinates $x=ry$ with $|y|=1$ to
find
$$
L^*L(P)
\ =\ 
|S^{d-1}|\
\int_0^R\
\int_{|y|=1}
r^{d-1}\
e^{-iry\xi}\
\phi_{d,k}( r )\
r^k\ P(y)
\
d\sigma(y)\, dr\,.
$$
Formula
\eqref{eq:FThomog3} shows that
$$
\int_{|y|=1}
e^{-iry\xi}\
P(y)
\
d\sigma(y)\
\ =\
\phi_{d,k}(r)\
P(-r\,\xi)
\ =\
(-r)^k\, 
\phi_{d,k}(r)\,
P(\xi)\,.
$$
Therefore
$$
L^*L(P)
\ =\ 
\int_0^R
(-1)^k\,
|S^{d-1}|\
r^{d+2k-1}\
\phi_{d,k}(r)^2\
dr\ P(\xi)\,.
$$
This proves that the harmonic polynomial are
eigenvectors with eigenvalue depending only 
on $d$ and $k$.

The formula for the eigenvalue follows on noting
that the eigenvalue is equal to the square of 
the $L^2(B_R(0))$ norm of 
\begin{equation}
\label{eq:defF}
F(x):=
\int_{|\eta|=1}
\
e^{ix\eta}
\
P(\eta)
\ 
d\sigma\,,
\qquad
\int_{|\eta|=1}
\big|P(\eta)\big|^2
\ 
d\sigma =1.
\end{equation}
Using polar coordinates and \eqref{eq:FThomog3}  yields,
\begin{equation}
\label{eq:concentration}
\begin{aligned}
\big\|
F\big\|_{B_R(0)}^2
\ &=\ 
(2\pi)^{d/2}\,
|S^{d-1}|
\ 
\int_0^R
r^{-(d-2)}
\ 
\big[J_{(d+2k-2)/2}(r)\big]^2
\ r^{d-1}
\ dr
\cr
\ &=\ 
(2\pi)^{d/2}\,
|S^{d-1}|
\ 
\int_0^R
r
\ 
\big[
J_{(d+2k-2)/2}(r)
\big]^2
\ dr
\,.
\end{aligned}
\end{equation}
This proves \eqref{eq:Lambdadef}.
\qed

\vskip.2cm

The spectral decomposition of $\bfL$ is nearly identical
to that of $L$.  The next result is elementary.

\begin{corollary} 
\label{cor:bfLspec} The eigenvalues of $\bfL^*\bfL$ are 
the same as the eigenvalues of $L^*L$.  The eigenspaces
consists of vector valued functions each of whose components
belongs to the corresponding eigenspace of $L^*L$.
\end{corollary}

\subsection{Some eigenfunctions and eigenvalues of $\bfL_M^*\bfL^{}_M$.}

The situation for $\bfL^*_M\bfL^{}_{M}$ is more subtle.  
Our first two results show that there are eigenfunctions
intimately related to the eigenvalus $\Lambda_{d,1}(R)$
and $\Lambda_{d,0}(R)$.

\begin{theorem}  
\label{thm:H1} The $d$ dimensional space of functions
$\bfe(\xi):=\zeta\wedge\xi$ with $\zeta\in \CC^d\setminus 0$
are eigenfunctions of $\bfL_M^*\bfL^{}_M$.   The eigenvalue is
$\Lambda_{d,1}(R)$.
\end{theorem}

\begin{remark}  These $\bfe(\xi)$ are the vector valued
spherical harmonics of degree 1 that belong to $\bfH$, that is,
that satisfy $\xi\cdot\bfe(\xi)=0$.
\end{remark}

{\bf Proof.}  Follows from 
$\bfL \bfe=\Lambda_{d,1}(R)\,\bfe$
and $\bfe\in \bfH$.
\qed

\vskip.2cm

Though the constant functions which are eigenvectors
of $\bfL$ do not belong to $\bfH$,
their projection on
$\bfH$ yield eigenvectors of $\bfL_M^*\bfL^{}_M$.

\begin{theorem}
\label{thm:LMspec}  The $d$ dimensional space
consisting of scalar multiples of rotates of 
$\ell(\xi)$ consists of eigenfunctions of $\bfL_M^*\bfL^{}_M$
  with eigenvalue equal to 
  \begin{equation}
  \label{eq:LMevalue}
  \big(\Lambda_{d,0}(R)-\Lambda_{d,2}(R)\big)\
  \frac{d-1}{d}
  \,.
 \end{equation}
\end{theorem}

  {\bf Proof of theorem.}  Use \eqref{eq:ellexp}.
Since the spherical harmonics are eigenfunctions of
$\bfL^*\bfL$ one has, suppressing the $R$ dependence of $\Lambda$,
\begin{equation*}
\label{LstarLell}
\bfL^*\bfL\,\ell
\ =\
 \Big(
   \Lambda_{d,0}\frac{d-1}{d} -
   \Lambda_{d,2}\sum_{j=2}^d
   \frac{\xi_1^2-\xi_j^2}{d}
\,,\,
-\Lambda_{d,2}\,\xi_1\xi_2
\,,\,
\cdots
\,,\,
-\Lambda_{d,2}\,\xi_1\xi_d
\Big).
\end{equation*}

Multiply \eqref{eq:ellexp} by $\Lambda_{d,2}$ to find
on $\xi_1^2+\cdots +\xi_d^2=1$,
$$  
   \Lambda_{d,2}\,
   \ell\ =\
    \bigg(
     \Lambda_{d,2}\,\frac{d-1}{d} -
    \Lambda_{d,2}\, \sum_{j=2}^d
   \frac{\xi_1^2-\xi_j^2}{d}
\,,\,
-  \Lambda_{d,2}\,\xi_1\xi_2
\,,\,
\cdots
\,,\,
-  \Lambda_{d,2}\,\xi_1\xi_d
\bigg)
\,.
$$
Subtract
from the preceding identity to find,
$$
 \bfL^*\,\bfL
\,
   \ell
\ =\  
\Big( \Lambda_{d,0}-\Lambda_{d,2}\Big)\frac{d-1}{d} 
  \big(1,0,\dots,0\big)
  \,.
  $$
  Projecting  perpendicular to $\xi$ 
  using $\Pi\,(1,0,\cdots,0)= \ell$ yields
    $$
  \Pi\,
  \bfL^*\,\bfL\,\Pi\,\ell
  \ =\
  \Pi\,
  \bfL^*\bfL\, \ell
  \ =\ 
 \Big( \Lambda_{d,0}-\Lambda_{d,2}\Big)
   \frac{d-1}{d} 
  \
  \ell
    \,.
  $$
  This proves that $\ell$ is an eigenfunction of 
  $(\bfL\Pi)^*(\bfL\Pi)$
  with eigenvalue $\big(\Lambda_{d,0}-\Lambda_{d,2}\big)(d-1)/d$.
  
  Remark \ref{rmk:LL}.iii  shows that it is an eigenfunction
  of $\bfL_M^*\bfL^{}_M$ with the same eigenvalue.

  By rotation invariance the same is true of all scalar multiples of
  rotates of $\ell$.  They form a $d$ dimensional vector space.
    spanned by the projections tangent to the unit sphere of the
  unit vectors along the coordinate axes.
  \qed

\vskip.2cm

As in Theorem \ref{thm:H1} 
 if one defines $\bfH_k$ to consist of spherical
harmonics of degree $k$ that belong to $\bfH$, then 
$\bfH_k$
are orthogonal eigenspaces of $\bfL^*_M\bfL^{}_M$
with eigenvalue $\Lambda_{d,k}(R)$.

\begin{example}
In $\RR^2$ the homogeneous $\RR^2$ valued
polynomials of degree two
whose radial components vanish 
are spanned by $(-x_1x_2, x_1^2)$
and $(x_2^2, -x_1x_2)$.  
There are no harmonic functions is their span
proving that 
when $d=2$, $\bfH_2=0$.  It is clear that
$\bfH_0=0$.

Though there are a substantial number  of eigenvectors
of $\bfL_M^*\bfL^{}_M$ accounted for by the $\bfH_k$ they are
far from the whole story.
\end{example}

\section{Spectral asymptotics.}

\subsection{Behavior of  the $ \Lambda_{d,k}(R)$.}

\begin{proposition}
\label{prop:Lambdaasym}
  {\bf i.} As $R\to 0$, 
$\ \ \Lambda_{d,k}(R)=O(R^{d+2k})$.

{\bf ii.}  As $R\to 0$, $\ \ \Lambda_{d,0}(R) =
| S^{d-1}|\,|B_R(0)|(1+O(R))$.

{\bf iii.}  $\lim_{k\to \infty} \Lambda_{d,k}(R)\, =\, 0\ $
uniformly on compact sets of $R$.
\end{proposition}

{\bf Proof.}  {\bf i.}  Formula \eqref{eq:Jdef} shows that
$J_k(t)=O(t^k)$ as $t\to 0$.  Assertion {\bf i}
then follows from \eqref{eq:Lambdadef}.

{\bf ii.}  By definition, $\Lambda_{d,0}(R)$ is the square
of the
$L^2(B_R)$ norm of $\int_{| \xi | =1} e^{ix\xi}\,f(\xi)\,d\sigma$
for $f$ a constant function of norm 1.  Take
$f=|S^{d-1}|^{-1/2}$.

For $R$ small $e^{ix\xi} = 1 + O(R)$ so
$$
\int_{| \xi | =1} e^{ix\xi}\,f(\xi)\,d\sigma
=
\int_{| \xi | =1} (1+O(R))\, |S^{d-1}|^{-1/2}\,d\sigma
=
|S^{d-1}|^{1/2}\big(1+O(R)\big)\,.
$$
Squaring and integrating over $B_R$ proves {\bf ii}.

{\bf iii.}  As $k\to \infty$,
the prefactors
in the formula for $J_k$ tend to zero uniformly since
$\Gamma([(2k+1)/2]$ dominates.

The integral in the definition  tends to zero 
uniformly by Lebesgue's Dominated Convergence Theorem.  
\qed

\subsection{Small $R$ asymptotics of the largest eigenvalues.}
\label{sec:smallR}

Each of the operators $L,\bfL,$ and $\bfL_M$ has
integral kernel $e^{ix\xi}$.  They differ in the Hilbert
space on which they act.  For $x$ small,
$e^{ix\xi}\approx 1$ showing that for $R$ small
the three operators are approximated by 
$L_0, \bfL_0,$ and $\bfL_{0,M}$ respectively.  
We know
 the
exact spectral decomposition of the
approximating operators.
Each has exactly one nonzero eigenvalue.
In performing the approximation some care must
be exercised since the operator to be approximated
has norm $O(R^{d/2})$ tending to zero as $R\to 0$.

\begin{proposition}
\label{prop:approxOp}
{\bf i.}  Each of the operators $L$, $\bfL$, and $\bfL_M$
has norm no larger than
$(|B_R(0)|\, |S^{d-1}|)^{1/2}$.

{\bf ii}  Each of the differences $L-L_0$,
$\bfL-\bfL_0$, and $\bfL_M-\bfL_{0,M}$ has norm
no larger than 
\begin{equation}
\label{eq:upbound}
\frac{
|S^{d-1}|\,
R^{(d+2)/2}
}
{
(d+2)^{1/2}
}
\,.
\end{equation}
\end{proposition}

{\bf Proof.}  {\bf i.}  
Treat the case of $\bfL$.  For $\| \bfe\| =1$, the Cauchy-Schwartz
inequality implies that for each $x$, $\|\bfL \bfe(x)\|_{\CC^d}^2
\le |S^{d-1}|$.  Integrating over the ball of radius $R$
proves {\bf i}.

{\bf ii}.  The Cauchy-Schwarz inequality estimates the difference by
$$
\begin{aligned}
\big|
(\bfL-\bfL_0)\bfe \big|
&=
\Big|
\int_{|\xi|=1}
\big(
e^{ix\xi}-1
\big)
\,
\bfe(\xi)
\,
d\sigma
\Big|
\le
\int_{|\xi|=1}
 |x| \, | \bfe |\, d\sigma
 \cr
&\le
|x| \, |S^{d-1}|^{1/2}\, \| \bfe \|_{L^2(S^{d-1})}
\,.
\end{aligned}
$$
For $ \bfe $ of norm one this yields
$$
\begin{aligned}
\|
(\bfL -\bfL_0) \bfe
&
\|_{L^2(B_R(0))}^2
\ \le\
|S^{d-1}|\,
\int_{B_R}
|x|^2\,dx\cr
\ &=\ 
|S^{d-1}|\,
\int_{|\omega|=1} 
\int_0^R
r^2\,r^{d-1}\, dr\, d\sigma(\omega)
\ =\
|S^{d-1}|^2\
\frac{R^{d+2}}{d+2}\,,
\end{aligned}
$$
completing the proof.
\qed

\begin{theorem}
\label{thm:smallR1}
{\bf i.}  For each $d$ there is an $R_L(d)>0$ so that
for $0\le R<R_L(d)$ the eigenvalue
$\Lambda_{d,0}(R)$
is the largest eigenvalue of $L^*L$.  
It has multiplicity one.  The eigenfunctions are
constants.  

{\bf ii.}  For $0\le R<R_L(d)$   the eigenvalue
$\Lambda_{d,0}(R)$
is the largest eigenvalue of $\bfL^*\bfL$.  
It has multiplicity d.   The eigenfunctions are
constant  vectors.

{\bf iii.}  For each $d$ there is an $R_M(d)>0$
so that
for $0\le R<R_M(d)$ the eigenvalue
 \begin{equation}
 \label{eq:evalbehave}
  \big(\Lambda_{d,0}(R)-\Lambda_{d,2}(R)\big)
  \frac{d-1}{d}
  \end{equation}
is the largest eigenvalue of $\bfL_M^*\bfL^{}_M$.  
It has multiplicity $d$.  The eigenfunctions are
rotates of constant multiples of $\ell$.

{\bf iv.}  In all three cases, the other eigenvalues are 
$
O(
R^{
d+1
}
)
$.
\end{theorem}

{\bf Proof.}  We prove {\bf iii} and {\bf iv}
for the operator $\bfL_M^*\bfL^{}_M$.
  Proposition
\ref{prop:approxOp} implies that 
the compact self adjoint operators
$\bfL_M^*\bfL^{}_M$ and 
$\bfL_{0,M}^*\bfL^{}_{0,M}$ 
differ by  $O(R^{d+1}
)$ in norm.

Part {\bf iii} of 
Theorem \ref{thm:Lzerospec} 
shows that the spectrum of 
$\bfL_{0,M}^*\bfL^{}_{0,M}$
contains one positive eigenvalue,
$\lambda_+:=|B_R(0)|\,|S^{d-1}|(d-1)/d$.
The factor $|B_R(0)|$ arises because
$\bfL_{0,M}$ in the present context is viewed
as an operator with values in the functions
on $B_R(0)\subset\RR^d$.  
The eigenfunctions are scalar multiplies of rotates
of $\ell$.  The rest of the spectrum
is 
the eigenvalue $0$.

It follows that the spectrum of 
$\bfL_{M}^*\bfL^{}_{M}$ lies in the 
union of disks of radius $O(R^{
d+1}
)$
centered at zero and
$\lambda_+$.
For $R$ small these disks are disjoint
and the eigenspace associated to
the disk about 
$\lambda_+$ 
has dimension $d$.

Theorem \ref{thm:LMspec} shows that
the eigenfunctions of $\bfL_{0,M}^*\bfL^{}_{0,M}$
with eigenvalue $\lambda_+$
are
eigenfunctions of $\bfL_{M}^*\bfL^{}_{M}$.
The eigenvalue is given by 
\eqref{eq:LMevalue}.

It follows that  for $R$ small,  the scalar multiples of rotates
of $\ell$ is an eigenspace of $\bfL_M^*\bfL^{}_M$
of dimension $d$ and eigenvalue in the disk
about
$\lambda_+$.
This completes the proof of {\bf iii}.

The fact that the other eigenvalues lie in a disk
of radius $O(R^{d+1})$ centered at the origin
proves {\bf iv} .

The proofs for the operators $L$ and $\bfL$ are
similar.
\qed

\section{Largest eigenvalues for $\bf R\le \pi/2$.}

\subsection{Monotonicity of $\bf \Lambda_{d,k}(R)$ in $\bf k$.}

Recall that the wavelength is equal to $2\pi$.

\begin{theorem}
\label{thm:smallR}
For $0\le R\le \pi/2$,  $\Lambda_{d,k}(R)$ is strictly
monotonically
decreasing in $k=0,1,2, \dots$.
In particular the largest eigenvalue of $L^*L$
and $\bfL^*\bfL$ is $\Lambda_{d,0}(R)$.
The corresponding
eigenfunctions are constant scalar
and constant vector functions respectively.
\end{theorem}

{\bf Proof.}  Write,
\begin{equation}
\label{eq:Jint}
J_k(r) \ =\ 
\frac{2\,(r/2)^k}{\Gamma[(2k+1)/2]\ \Gamma(1/2)}\
\int_{0}^1
\cos(rs)\
(1-s^2)^{(2k-1)/2}\
ds\,.
\end{equation}

For $0\le r\le \pi/2$ the cosine factor in the integral is positive.
Since $(1-s^2)^{(2k-1)/2}$ 
is decreasing in $k$ for $s\in [0,1]$,
the integral is
decreasing  in $k$.

$\Gamma[(2k+1)/2]$ is increasing in $k$.  
Since $r\le 2$, $(r/2)^k$ decreases with $k$.
The proof is complete.
\qed

\begin{remark}  The first figure in 
\S 
\ref{sec:scalarsim}
shows that
the graphs of the functions $\Lambda_{3,0}(R)$
and 
and $\Lambda_{3,1}(R)$ cross close to $R=\pi$. The proof
shows that this cannot happen for $R\le \pi/2$.  
\end{remark}

\subsection{Largest eigenvalues of $\bfL_M^*\bfL^{}_M$ when
$\bf d=3$, $\bf R\le \pi/2$.}

\begin{theorem}
\label{thm:Maxlargest}
When $d=3$ and $R\le \pi/2$ the strictly largest eigenvalue of
$\bfL_M^*\bfL_M$ is $(\Lambda_{3,0}(R)-\Lambda_{3,2}(R))2/3$.
The eigenfunctions are the scalar multiples of 
rotates of $\ell$,  and, 
$\Lambda_{3,1}(R)$ is the next largest
eigenvalue.
\end{theorem}

The proof uses the following criterion valid for all $d,R$.
For ease of reading, the $R$ dependence of $\Lambda_{d,k}(R)$
is often suppressed.

\begin{theorem}
\label{thm:evcond}
{\bf i.}  The eigenvalue 
$(\Lambda_{d,0}-\Lambda_{d,2})(d-1)/d$ of
$\bfL_M^*\bfL^{}_M$ 
is strictly larger than all others only if 
\begin{equation}
\label{eq:evcond}
(\Lambda_{d,0}-\Lambda_{d,2})(d-1)/d
\ >\
\Lambda_{d,1}
\,.
\end{equation}

{\bf ii.}  If in addition to \eqref{eq:evcond},
the two largest evalues of  $L^*L$
are 
$\Lambda_{d,0}$ and $\Lambda_{d,1}$,
then
the eigenvalue
$(\Lambda_{d,0}-\Lambda_{d,2})(d-1)/d$
  of
$\bfL_M^*\bfL^{}_M$  
is strictly larger than the others.  
The
eigenfunctions are the scalar multiples of 
rotates of $\ell$, and, $\Lambda_{d,1}$ is the next largest
eigenvalue of $\bfL_M^*\bfL^{}_M$.
\end{theorem}

\begin{remark}  Equation \eqref{eq:evcond} implies
that $\Lambda_{d,0}>\Lambda_{d,1}$.  The
additional condition in {\bf ii} is that for all
$k\ge 1$, $\Lambda_{d,1}\ge  \Lambda_{d,k}$.
\end{remark}

{\bf Proof of Theorem \ref{thm:evcond}.}  {\bf i.}  Since   $\Lambda_{d,1}$
is also an eigenvalue
of 
$\bfL_M^*\bfL^{}_M$, necessity is clear.

{\bf ii.}  Under these hypotheses
Theorem \ref{thm:lambdaofR}
shows that $\Lambda_{d,0}$ is the largest
eigenvalue of $L^*L$ with one dimensional
eigenspace consisting of constant functions.
Corollary \ref{cor:bfLspec} shows that
$\bfL^*\bfL$ has the same largest eigenvalue
with $d$ dimensional eigenspace consisting
of $\CC^d$ valued constant functions.
The next largest eigenvalue of $\bfL^*\bfL$ is
$\Lambda_{d,1}$.  In particular,
$\bfL^*\bfL$ has exactly $d$ eigenvalues
counting multiplicity that are greater than
$\Lambda_{d,1}$.

Since $\bfL_M$ is the restriction of $\bfL$
to a closed subspace,
the minmax principal implies that 
$\bfL_M^*\bfL^{}_M$ 
has at most
$d$ 
eigenvalues
counting multiplicity that are greater than
$\Lambda_{d,1}$.

Theorem \ref{thm:LMspec} provides a 
$d$ dimensional eigenspace with 
eigenvalue given by the left hand side
of \eqref{eq:evcond} and therefore
greater than $\Lambda_{d,1}$.  In particular
there are exactly $d$ eigenvalues counting
multiplicity that are greater than 
$\Lambda_{d,1}$.

Theorem \ref{thm:H1} shows that $\Lambda_{d,1}$
is an eigenfunction of $\bfL_M^*\bfL^{}_M$ so it must
be the next largest.
\qed

\begin{example}  Parts {\bf i} and {\bf ii}
of Proposition \ref{prop:Lambdaasym}
show that the sufficient condition is satisfied 
for small $R$.  This gives a second
proof that for small $R$, $\ell$ is an extreme eigenfunction
for $\bfL_M^*\bfL_M$.  The first proof is part
{\bf iii} of Theorem \ref{thm:smallR1}
\end{example}

{\bf   Proof of Theorem \ref{thm:Maxlargest}.}    Verify the sufficient
condition of Theorem \ref{thm:evcond}.ii.

Since $R\le \pi/2$, $\Lambda_{3,k}$ are strictly
decreasing in $k$.  Therefore $\Lambda_{3,0}$
and $\Lambda_{3,1}$ are the two largest eigenvalues
of $L^*L$.

It remains to verify
\eqref{eq:evcond}.
Use formulas
\eqref{eq:Lambdadef}
and
\eqref{eq:Jint}.
Formula \eqref{eq:Lambdadef} with $d=3$ and $k=0,1,2$
involves $J_{1/2}, J_{3/2}, J_{5/2}$.

Since the integral in \eqref{eq:Jint} is decreasing in $k$
it follows that
$$
\frac{J_{k+1} (r)}
{ J_k (r)}
\ \le\
\frac{ r } { 2}\
\frac{ \Gamma( (2k+1)/2 ) }
{\Gamma( (2k+1)/2\, +1 ) }
\,.
$$
The functional equation $\Gamma(n+1)=(n+1)\Gamma(n)$
yields
$$
\frac{J_{k+1} (r) }
{ J_k (r) }
\ <\
\frac{ r } { 2}\
\frac{1}
{(2k+3)/2}
\ =\ 
\frac{r}
{2k+3}
\,.
$$
Therefore,
$$
\frac{J_{3/2} (r) }
{J_{1/2} (r) }
\ <\
\frac{r}{4},
\qquad
{\rm and},
\qquad
\frac{J_{5/2} (r) }
{J_{3/2} (r) }
\ \le\
\frac{r}{6}
\,.
$$

Injecting these estimates in \eqref{eq:Lambdadef}
yields
\begin{equation}
\label{eq:Lambdadecreaase}
\frac{ \Lambda_{3,1}(R) }
{ \Lambda_{3,0}(R) }
\ <\
 \frac{R^2}{4^2}
\,,
\qquad
{\rm and},
\qquad
\frac{ \Lambda_{3,2}(R) }
{ \Lambda_{3,1}(R) }
\ <\
 \frac{R^2}{6^2}
\,.
\end{equation}
Therefore,
$$
\Lambda_{3,1}
\ \le\
 \frac{R^2}{4^2}\,
 \Lambda_{3,0},
 \quad
 {\rm and},
 \quad
 \Lambda_{3,2}
 \ \le\
  \frac{R^2}{6^2}\,
 \Lambda_{3,1}
\ \le\
 \frac{R^2}{4^2}\,
  \frac{R^2}{6^2}\,
 \Lambda_{3,0}\,,
 \quad
 {\rm so}\,,
 $$
$$
\big(\Lambda_{3,0} -
\Lambda_{3,2}\big)\frac{2}{3} -\Lambda_{3,1}
\ >\
\Lambda_{3,0}\bigg(
\frac{2}{3}
\ -\
\frac{2}{3}\,
\frac{R^4}{4^2 \, 6^2}
\ -\
\frac{R^2}{4^2}
\bigg)
\ :=\ 
\Lambda_{3,0}\,h(R)\,.
$$
The polynomial $h(R)$ is equal to 
$2/3$ when $R=0$ and decreases as $R$
increases.  To verify \eqref{eq:evcond} it 
suffices to
show that $h(\pi/2)>0$.
Since $2>\pi/2$, $h(\pi/2)>h(2) = 43/108>0$.
\qed

\section{Numerical simulations to determine largest eigenvalues.}

Recall that the wavelength
is equal to $2\pi$.  In this section the dimension  $d=3$.

Theorem \ref{thm:lambdaofR},
Corollary \ref{cor:bfLspec}, and Theorem
\ref{thm:Maxlargest} allow one in favorable cases
to find the largest eignevalues of $L^*L$, $\bfL^*\bfL$,
and 
$\bfL_M^*\bfL^{}_M$, by evaluating the intergrals
defining $\Lambda_{3,k}(R)$ for $k=0,1,2,\dots$.
These quantities decrease rapidly
with $k$ so to compute  the largest ones requires little
work.

\subsection{Simulations for scalar waves.}
\label{sec:scalarsim}

For scalar waves the eigenvalues are
exactly 
the $\Lambda_{3,k}(R)$.
For $R$ small, they are monotone in $k$
so the optimal focusing is for $k=0$.
Our first simulation (performed with the aid of
Matlab) computes approximately the integrals
defining $\Lambda_{3,k}(R)$ for $R\le 2\pi$
and $k=0,1,2,3$.  The resulting graphs
are in the figure on the left.  The horizontal axis
is $R$ and 
on the vertical axis is plotted
the integral on the right hand side 
of 
\eqref{eq:Lambdadef}, that is,
$$
\frac
{\Lambda_{3,k}(R)}
{
(2\pi)^{3/2}|S^2|
}
\ =\
\frac
{\Lambda_{3,k}(R)}
{ 2^{7/2}\ \pi^{5/2} }
\,,
\qquad
k=0,1,2,3
\,.
$$

\begin{figure}
\begin{center}
 \includegraphics[height=45mm]{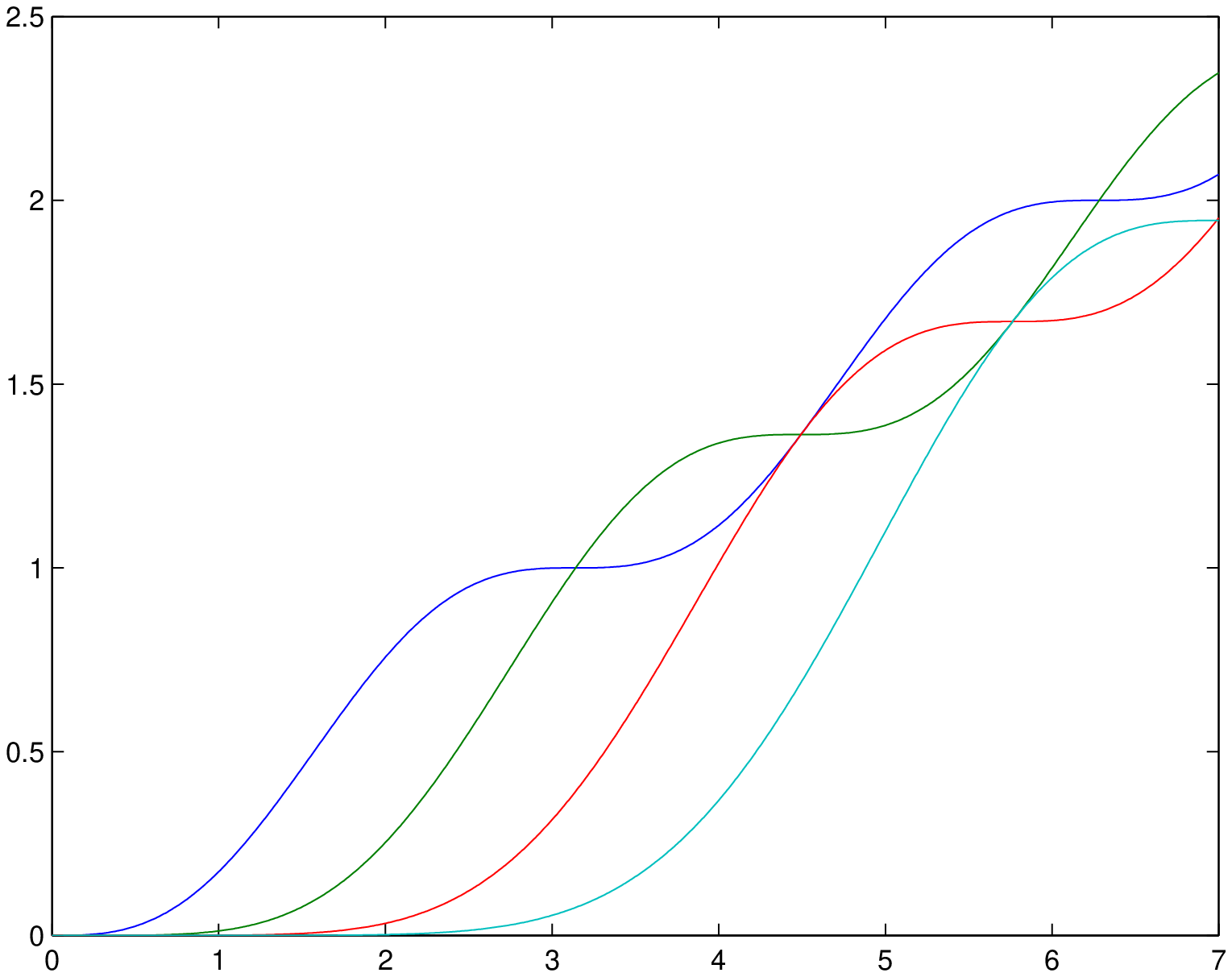}
 \hskip.3cm
\includegraphics[height=45mm]{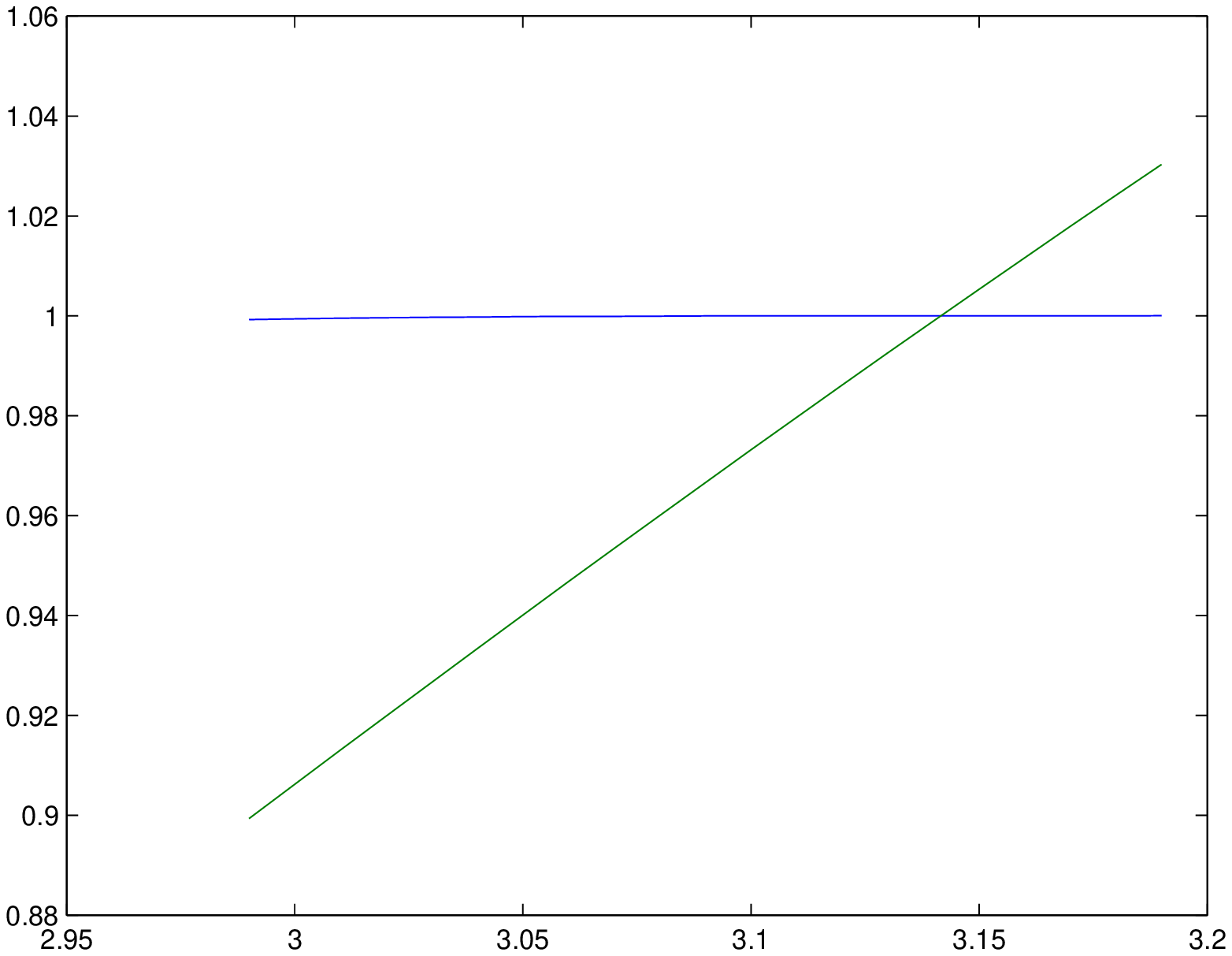}
\end{center}
\end{figure}

The four curves correspond to the four values of $k$.
The graph with the leftmost
hump is $\Lambda_{3,0}(R)$.  The graph  with the hump
second from the left is $\Lambda_{3,1}(R)$ and so on.
The conclusion is that $\Lambda_{3,0}(R)$  crosses
transversaly the graph
$\Lambda_{3,1}(R)$ just to the right of $R=3$.  At that
point, $\Lambda_{3,1}(R)$ becomes the largest.
 On the right
is a zoom showing  
that the crossing is suspiciously close to $R=\pi$.

The graphs of the $\Lambda_{3,k}(R)$ are a little misleading since
it is not the total energy but the energy density that
is of interest.  The next figure plots as a function of 
$R$ $\Lambda_{3,0}/(2^{7/2}\pi^{5/2}|B_R(0)|)$.
The small gap near $R=0$ is because the division by
$|B_R(0)|$ is a sensitive operation and leads to numerical
errors in that range.

\begin{figure}
\begin{center}
 \includegraphics[height=45mm]{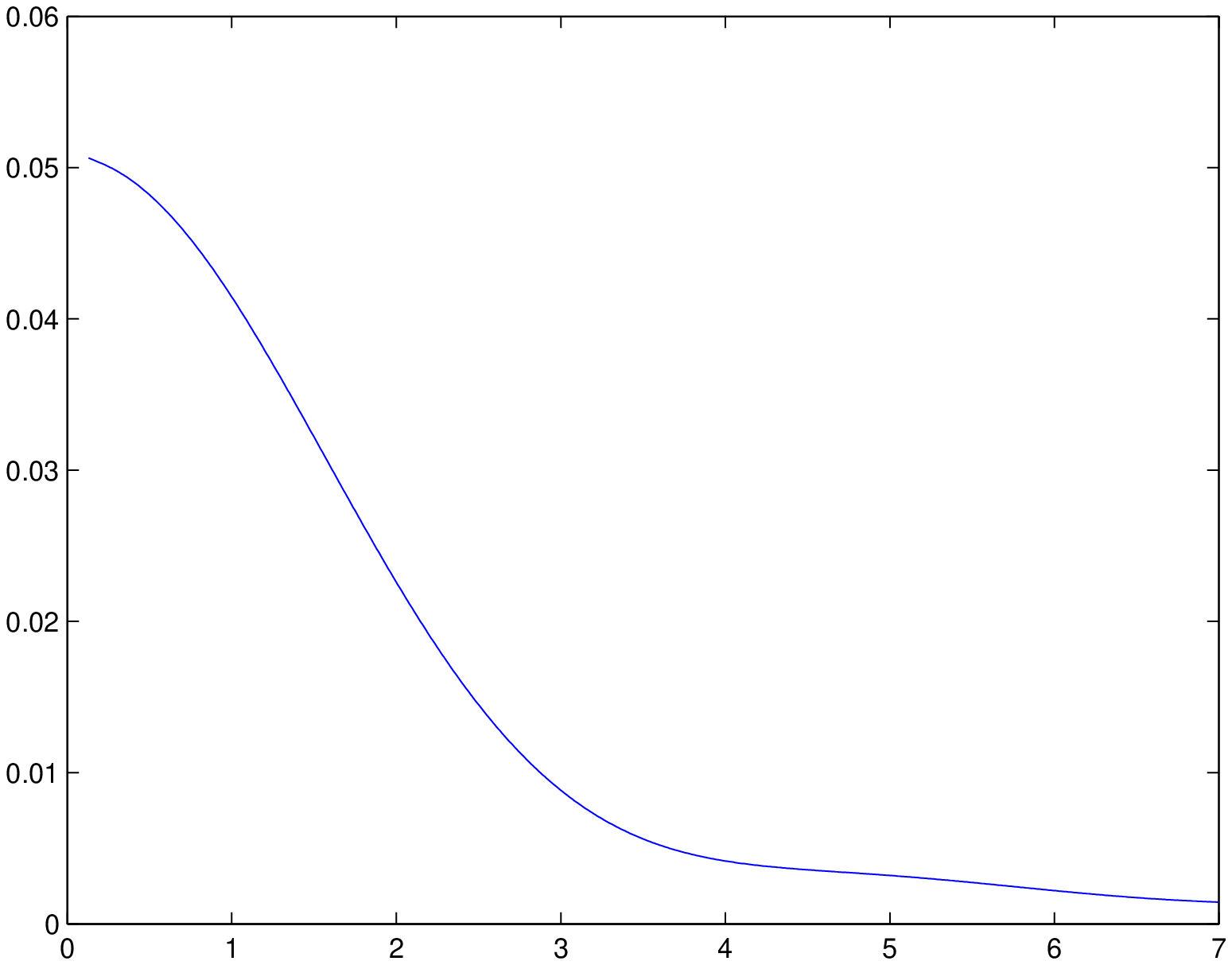}
\end{center}
\end{figure}

The energy density is greatest for balls with radius
close to $R=0$.   The density drops to half its maximum
value at about $R=2$ which is about 1/3 of the wavelength.

 \subsection{Simulations for electromagnetic waves.}

Using Theorem \ref{thm:Maxlargest}
one can investigate the analogous questions
for Maxwell's equations by manipulations of the
$\Lambda_{3,k}(R)$.

The simulations of the preceding subsection show that 
for $R\le 3$ one has $\Lambda_{3,0}(R)>\Lambda_{3,1}(R)$.
To show that the eigenvalue corresponding to $\ell(\xi)$
is the optimum it suffices to verify 
\eqref{eq:evcond}.
To do so one needs to verify the positivity of 
$$
2^{-7/2}\pi^{-5/2}
\bigg(\frac{2}{3}\,
\frac{
\Lambda_{3,0}(R)
}
{
|B_R(0) |
}
\ -\
\frac{2}{3}\,
\frac{
\Lambda_{3,2}(R)
}
{
|B_R(0) |
}
\ -\
\frac{
\Lambda_{3,1}(R)
}
{
|B_R(0) |
}
\bigg)
\,.
$$
This is a linear combination of quantities
computed in the 
preceding subsection.  Its graph is plotted on the left.   The graph
 crosses from positive to negative
near $R=2.5$.  The criterion is satisfied for all $R$
to the left of this crossing.

\begin{figure}
\begin{center}
  \includegraphics[height=45mm]{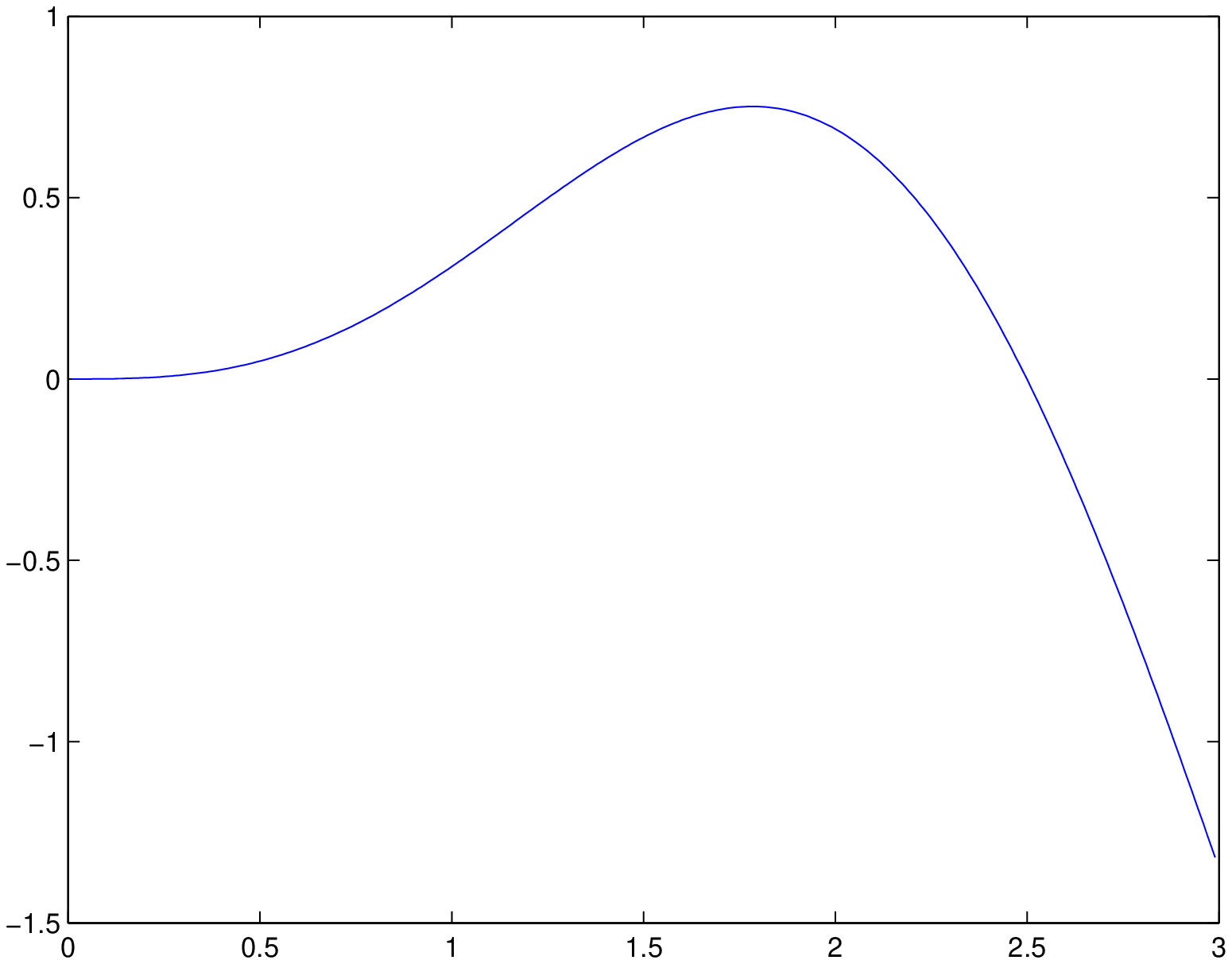}
  \hskip.3cm
\includegraphics[height=45mm]{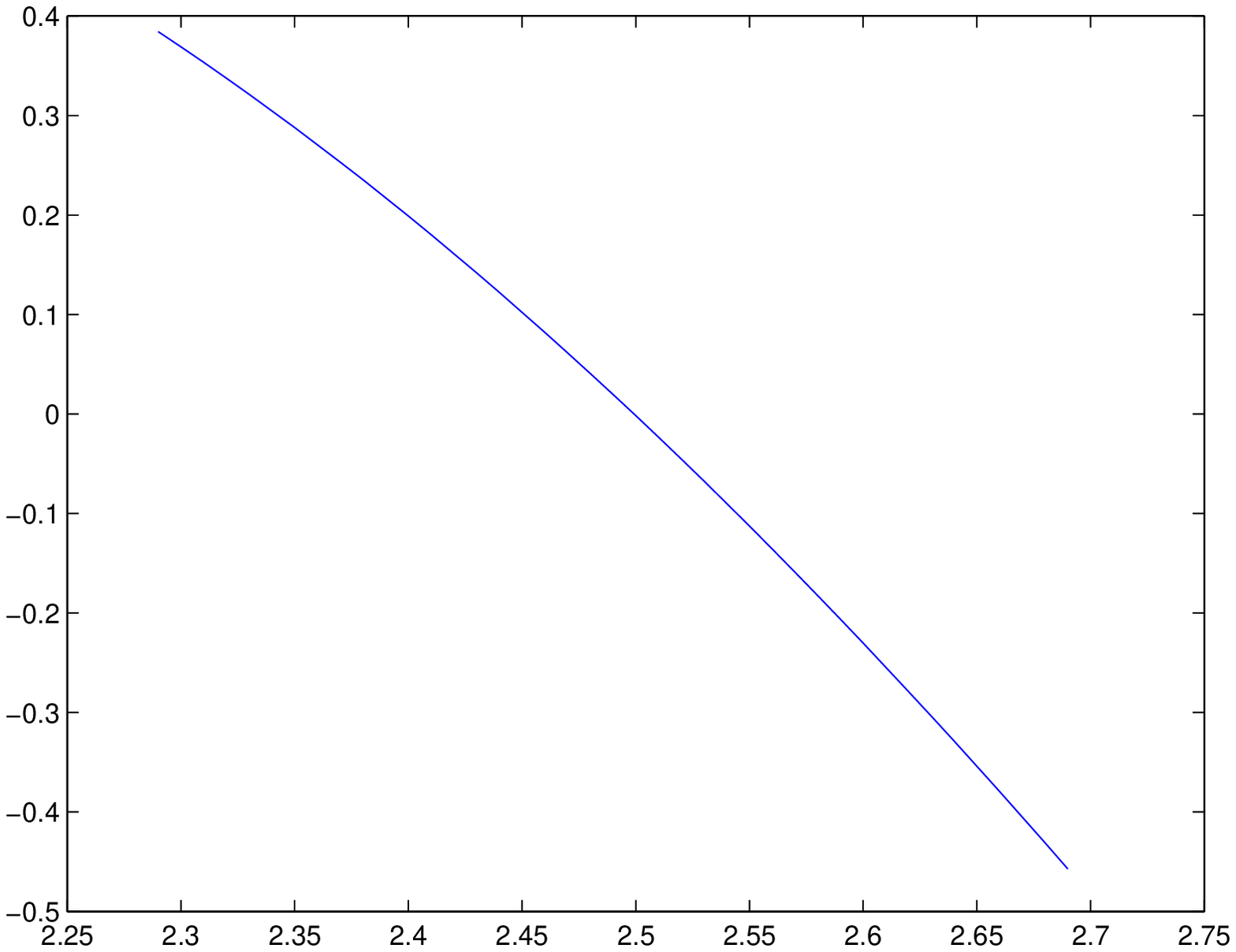}
\end{center}
\end{figure}

For $R<2.5$  the  energy density for the 
optimizing monochromatic 
electromagnetic fields associated with $\ell(\xi)$
is equal to
$$
\frac{2}{3}\,
\frac{
\Lambda_{3,0}(R)
}
{
|B_R(0) |
}
\ -\
\frac{2}{3}\,
\frac{
\Lambda_{3,2}(R)
}
{
|B_R(0) |
}
$$
Because of the factor $2/3$ it is smaller than the density
in the scalar case by that factor.  The subtraction in the formula
shows that the density drops off more rapidly in the electromagnetic
case than in the scalar case.
The graph of $2^{-7/2}\pi^{-5/2}$ times this quantity is
plotted next.

% \begin{figure}
\begin{center}
  \includegraphics[height=45mm]{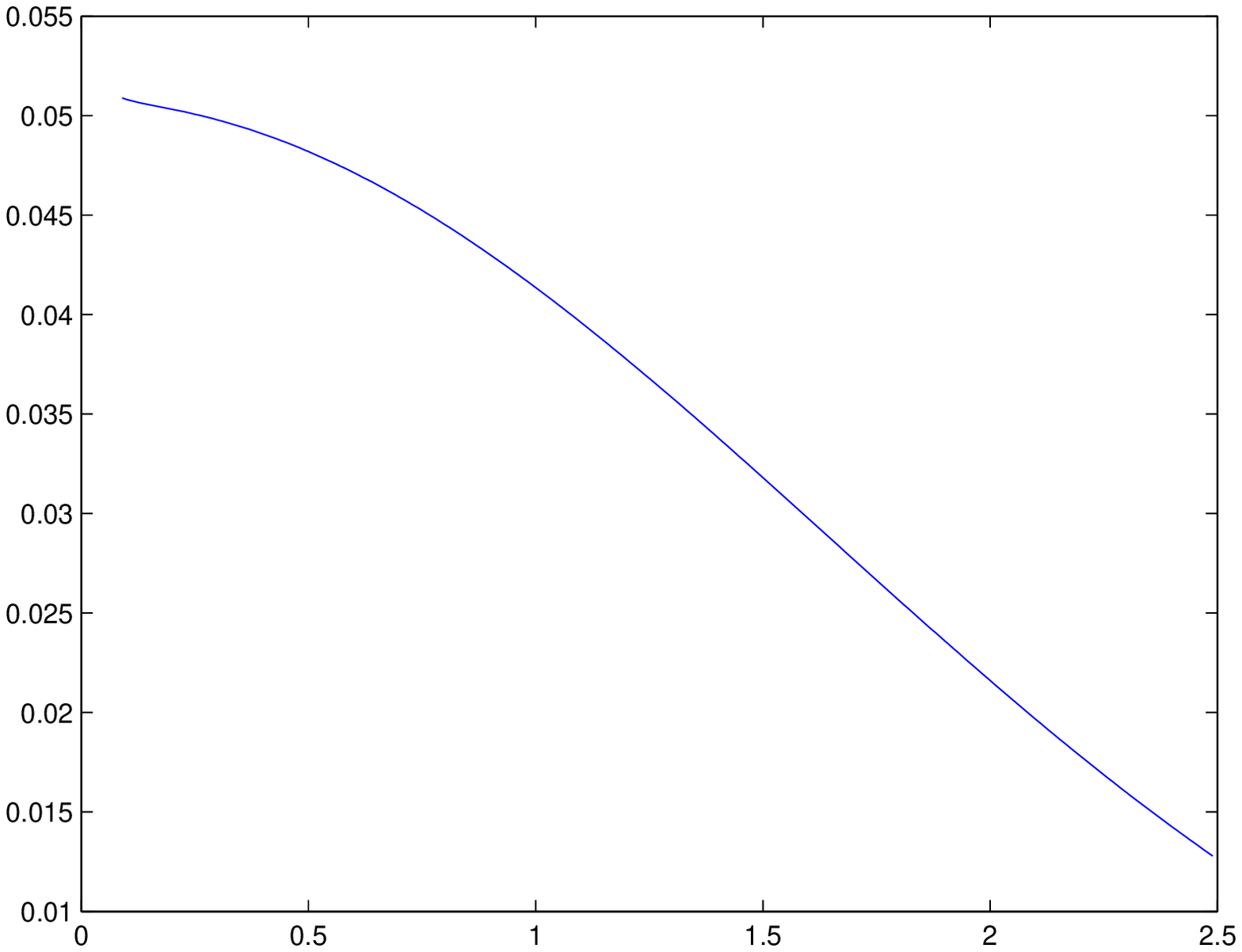}
\end{center}
% \end{figure}

As in the scalar case the maximal energy density occurs
on balls with radius near zero.  The intensity drops to one
half of this value to the left of $R=1.9$.
This is close to the 
corresponding value for scalar waves,
 about one third of a wavelength.

\end{document}